\documentclass[nonblindrev]{informs3_noinforms}
\usepackage{pdfsync}
\usepackage{enumerate, float, tikz}
\usepackage{srcltx}

\OneAndAHalfSpacedXI 

\usepackage{endnotes}
\let\footnote=\endnote

\def\prob {{\sf Pr}}

\newcommand{\cX}{{\mathfrak X}}

\newcommand{\cP}{{\mathfrak P}}

\newcommand{\cM}{{\mathfrak M}}
\newcommand{\cY}{{\mathfrak Y}}

\newcommand{\cf}{{\mathfrak f}}

\newtheorem{lemma}{Lemma}[section]
\newtheorem{proposition}{Proposition}[section]
\newtheorem{example}{Example}
\newtheorem{definition}{Definition}[section]

 \newcommand{\ii}{{\cal I}}

\newcommand{\I}{{\cal I}}

\newcommand{\half}{ \mbox{\small$\frac{1}{2}$}}

\def\val{{\rm val}}

\def\supp {{\rm supp}}

\def\sol {{\rm Sol}}

\def\bbr{{\Bbb{R}}} 
\def\bbe{{\Bbb{E}}} 

\newtheorem{theorem}{Theorem}
\newtheorem{observation}{Observation}

\renewcommand{\theequation}{\thesection.\arabic{equation}}

\usepackage[numbers]{natbib}
 \bibpunct[, ]{(}{)}{,}{a}{}{,}%
 \def\bibfont{\small}%
\ECRepeatTheorems

\EquationsNumberedThrough    

\begin{document}




\TITLE{Time (in)consistency of multistage distributionally robust inventory models with moment constraints}
\RUNTITLE{Time (in)consistency of multistage robust inventory models}


\ARTICLEAUTHORS{
\AUTHOR{\bf Linwei Xin}
\AFF{Booth School of Business, University of Chicago, Chicago, IL 60637}
\AUTHOR{\bf David A. Goldberg}
\AFF{Operations Research and Information Engineering, Cornell University, Ithaca, NY 14850}
\RUNAUTHOR{{Xin and Goldberg}}
}

\ABSTRACT{
\indent Recently, there has been a growing interest in developing inventory control policies which are robust to model misspecification.  One approach is to posit that nature selects a worst-case distribution for any relevant stochastic primitives from some pre-specified family.  Several communities have observed that a subtle phenomena known as time inconsistency can arise in this framework.  In particular, it becomes possible that a policy which is optimal at time zero (i.e. solution to the multistage-static formulation) may not be optimal for the associated optimization problem in which the decision-maker recomputes her policy at each point in time (i.e. solution to the distributionally robust dynamic programming formulation), which has implications for implementability.  If there exists a policy which is optimal for both formulations (w.p.1 under every joint distribution for demand belonging to the uncertainty set), we say that the policy is \emph{time consistent}, and the problem is \emph{weakly time consistent}.  If every optimal policy for the multistage-static formulation is time consistent, we say that the problem is \emph{strongly time consistent}.
\\\indent We study these phenomena in the context of managing an inventory over time, when only the mean, variance, and support are known for the demand at each stage.  We provide several illustrative examples showing that here the question of time consistency can be quite subtle.  We complement these observations by providing simple sufficient conditions for weak and strong time consistency.  We also relate our results to the well-studied notion of \emph{rectangularity} of a family of measures. Interestingly, our results show that time consistency may hold even when rectangularity does not.  Although a similar phenomena was previously identified by Shapiro for the setting in which only the mean and support of the demand are known, there the problem was always weakly time consistent, with both formulations having the same optimal value.  Here our model is rich enough to exhibit a variety of interesting behaviors, including lack of weak time consistency, strong time consistency even when both formulations have different optimal values, and non-existence of even a single optimal base-stock policy under the static formulation. 
}



\maketitle

\section{Introduction}\label{sec-introduction}

The news vendor problem, used to analyze the trade-offs associated with stocking an inventory, has its origin in a seminal paper by \cite{edg:88}.  In its classical formulation, the problem is stated as a minimization of the expected value of the relevant ordering, backorder, and holding costs. Such a formulation requires a complete specification of the probability distribution of the underlying demand process.  However, in applications knowledge of the exact distribution of the demand process is rarely available. This motivates the study of minimax type (i.e. distributionally robust) formulations, where minimization is performed with respect to a worst-case distribution from some family of potential distributions. In a pioneering paper \cite{Scarf} gave an elegant solution for the minimax news vendor problem when only the first and second order moments of the demand distribution are known.  His work has led to considerable follow-up work (cf. \cite{GM,moon1994distribution,gallego1998new,gallego2001minimax,Popescu,yue2006expected,GKR,PR,ChenSim,SS, hanasusanto2012distributionally,zhu2013newsvendor}).  For a more general overview of risk analysis for news vendor and inventory models we can refer, e.g., to \cite{ACS} and \cite{CRZ:11}.  We also note that a distributionally robust minimax approach is not the only way to model such uncertainty, and that there is a considerable literature on alternative approaches such as the robust optimization paradigm (cf. \cite{kasugai1961note,BGNV,BT06,aharon2009robust,BIP,carrizosaa2014robust,VMT14}) and Bayesian approach (cf. \cite{scarf1959bayes, scarf1960some, Lovejoy, LPU11, KSLS}).

In practice an inventory must often be managed over some time horizon, and the classical news vendor problem was naturally extended to the multistage setting, for which there is also a considerable literature (see, e.g., \cite{Zi} and the references therein).  Recently, distributionally robust variants of such multistage problems have begun to receive attention in the literature (cf. \cite{gallego2001minimax, ACS, choi2008risk, SS, S-12, KSLS}).  It has been observed that such multistage distributionally robust optimization problems can exhibit a subtle phenomenon known as time inconsistency.  Over the years various concepts of time consistency have been discussed in the economics  literature, in the context of rational decision making.  This can be traced back at least to the work of \cite{Strotz55} - for a more recent overview we refer the reader to the recent survey by \cite{EJT}, and the references therein.  Questions of time consistency have also attracted attention in the mathematical finance literature, in the context of assessing the risk and value of investments over time, and have played an important role in the associated theory of coherent risk measures (cf. \cite{Wan99,ADEHK,RS07,CK2009,Ru}).
These concepts have also been studied from the perspective of robust control across various academic communities (cf. \cite{HS01, iye:05,nil:05,GH2011,carp,WKR}).  Recently, these concepts have also begun to receive attention in the setting of inventory control (cf. \cite{CSLS,ChenSun,Yang13, Homem.16, Shapiro17}).

In this work, we will consider questions of time (in)consistency in the context of managing an inventory over time.  We will give a formal definition of time consistency, which is naturally suited to our framework,  in Section\ \ref{consistencysection}. At this point let us provide the following high-level intuition.  A multistage distributionally robust optimization problem can be viewed in two ways.  In one formulation, the policy maker is allowed to recompute her policy choice after each stage (we will refer to this as the distributionally robust dynamic programming (DP) formulation), thus taking prior realizations of demand into consideration when performing the relevant minimax calculations at later stages. In that  case it follows from known results that there exists a base-stock policy which  is optimal.  In the second formulation, the policy maker is not allowed to recompute her policy after each stage (we will refer to this as the {\em multistage-static formulation}), in which case far less is known.  If these two formulations have a common optimal policy, i.e. the policy maker would be content with the given policy whether or not she has the power to recompute after each stage (w.p.1 under every joint distribution for demand belonging to the uncertainty set), we say that the policy is \emph{time consistent}, and the problem is \emph{weakly time consistent}.  If every optimal policy for the multistage-static formulation is time consistent, i.e. it is impossible to devise a policy which is optimal at time zero yet suboptimal at a later time, we say that the problem is \emph{strongly time consistent}.  Such a property is desirable from a policy perspective, as it ensures that previously agreed upon policy decisions remain rational when the policy is actually implemented, possibly at a later time.

Within the optimization and inventory control communities, much of the work on time consistency restricts its discussion of optimal policies to the setting in which the family of distributions from which nature can select satisfies a certain factorization property called \emph{rectangularity}, which endows the associated minimax problem with a DP structure.  Outside of this setting, most of the literature focuses on discussing and demonstrating hardness of the underlying optimization problems (cf. \cite{iye:05,nil:05,WKR}).  We note that this is in spite of the fact that previous literature has discussed the importance and relevance of such non-decomposable formulations from a modeling perspective (cf. \cite{iye:05}).

\subsection{Our contributions}
In this paper, we depart from much of the past literature by seeking both negative \emph{and positive} results regarding time consistency when no such decomposition holds, i.e. the underlying family of distributions from which nature can select is non-rectangular.  Our work is in the spirit of \cite{GH2011}, in which a definition of (weak) time consistency similar to ours was analyzed in the context of rectangularity and dynamic consistency (a concept defined in \cite{ES}), albeit in a substantially different context motivated by questions in decision theory and artificial intelligence.  Our work can also be viewed as providing a more in-depth and inventory-focused study of several notions of time-consistency studied in \cite{Homem.16}.  In contrast to \cite{Homem.16} and several other works in which all concepts are explained through the language of risk measures, here we explain all relevant concepts purely in the language of (robust) newsvendor models with moment constrainst, a model popular in the operations management community, and hope that in doing so our work brings the concept of time-consistency to a broader audience.

We extend the work of \cite{Scarf} (and followup work of \cite{gallego2001minimax}) by considering the question of time consistency in multistage news vendor problems when the support and first two moments are known for the demand at each stage, and demand is stage-wise independent.  In addition to refining multiple definitions related to time-consistency, we provide several illustrative examples showing that here the question of time consistency can be quite subtle.  In particular: (i) the problem can fail to be weakly time consistent, (ii) the problem can be weakly but not strongly time consistent, and (iii) the problem can be strongly time consistent even if every associated optimal policy takes different values under the multistage-static and distributionally robust DP formulations.  We also prove that, although the distributionally robust DP formulation always has an optimal policy of base-stock form, there may be no such optimal policy for the multistage-static formulation.  We complement these observations by providing simple sufficient conditions for weak and strong time consistency.

Interestingly, in contrast to much of the related literature, our results show that time consistency may hold even when rectangularity does not.  This stands in contrast to the analysis of \cite{S-12} for the setting in which only the mean and support of the demand distribution are known, where the problem is always (weakly) time consistent, amenable to a simple DP solution, with both formulations having the same optimal value.  Likewise, in the setting in which only the support is known, both formulations reduce to the so-called adjustable robust formulation described in \cite{BGGN}, where again (weakly) time consistency always holds.  
\subsection{Outline of paper}
The structure of the rest of the paper is organized as follows.  In Section\ \ref{singlestagesec}, we review the single-stage classical and distributionally robust formulations and their properties, as well as Scarf's solution to the single-stage distributionally robust formulation and various generalizations.  In Section\ \ref{multisec}, we discuss the extension to the multi-stage setting, formally defining the multistage-static formulation, the relevant notions of time-consistency, and the distributionally robust DP formulation, and review the notion of rectangularity and its relation to our own formulations.  In Section\ \ref{consistencysection}, we prove our sufficient conditions for weak and strong time consistency, and present several illustrative examples showing that here the question of time consistency can be quite subtle.  In Section\ \ref{concsec}, we provide closing remarks and directions for future research.  We include a technical appendix in Section\ \ref{appsec}.

\section{Single-stage formulation}\label{singlestagesec}
In this section we review both the classical and distributionally robust single-stage formulation, including some relevant results of \cite{Scarf} and \cite{Natarajan}.
\subsection{Classical formulation}
\label{sec-stat}

Consider the following classical formulation of the news vendor problem:
\begin{equation}\label{nvc1}
\inf_{x\ge 0}\bbe [\Psi(x,D)],
\end{equation}
where
\begin{equation}\label{nvc1-a}
\Psi(x,d):=cx+b[d-x]_+ +h[x-d]_+,
\end{equation}
and $c,b,h$ are the ordering, backorder penalty, and holding costs, per unit, respectively. Unless stated otherwise  we assume  that $b>c> 0$ and  $h \ge  0$.
The expectation is taken with respect to the probability distribution of the demand $D$, which is modeled as a random variable   having nonnegative support.  It is well known that this problem has the closed form solution
$
\bar{x} = F^{-1}\left(\frac{b-c}{b+h}\right),
$
where $F(\cdot)$ is the cumulative distribution function (cdf) of the demand $D$, and $F^{-1}$ is its inverse.  Of course, it is assumed here that the probability distribution, i.e. the cdf $F$, is completely specified.
\subsection{Distributionally robust formulation}
\label{sec-disrob}
Suppose now that the probability distribution of the demand $D$ is not fully specified, but instead
assumed to be a member of a family of distributions $\cM$.  Then we  consider the following distributionally robust formulation:
\begin{equation}\label{nvc2}
\inf_{x\ge 0} \psi(x),
\end{equation}
where
\begin{equation}\label{eq-psi}
 \psi(x) :=\sup_{Q\in \cM}\bbe_Q [\Psi(x,D)],
\end{equation}
and the notation $\bbe_Q$ emphasizes that the expectation is taken with respect to the distribution $Q$ of the demand $D$.

We now introduce some additional notations to describe certain families of distributions.  For a probability measure (distribution) $Q$, we let $\supp(Q)$ denote the support of the measure, i.e. the smallest closed set $A\subseteq \bbr$ such that   $Q(A) = 1$.  With a slight abuse of notation, for a random variable $Z$, we also let $\supp(Z)$ denote the support of the associated probability measure.  For a given closed  (and possibly unbounded) subset $\ii\subseteq \bbr$, we let $\cP(\ii)$ denote the set of probability distributions $Q$ such that $\supp(Q)\subseteq \ii$.  Although we will be primarily interested in the setting that $\ii \subseteq \bbr_+$ (i.e. demand is nonnegative), it will sometimes be convenient for us to consider more general families of demand distributions. By $\delta_a$ we denote the probability measure of mass one at $a\in \bbr$.

In this paper, we will study families of distributions satisfying moment constraints of the form
\begin{equation}\label{nvc3}
\cM:=\big\{Q\in \cP(\ii):\bbe_Q[D]= \mu, \bbe_Q[D^2] = \mu^2 + \sigma^2\big\}.
\end{equation}
Unless stated otherwise, it will be assumed that $\cM$ is indeed of the form (\ref{nvc3}), and is nonempty.  We let $\alpha$ denote the left-endpoint of $\ii$ (or $-\infty$ if $\ii$ is unbounded from below), and let $\beta$ denote the right-endpoint of $\ii$ (or $+ \infty$ if $\ii$ is unbounded from above); i.e., $\ii = [\alpha,\beta]$.  Here we note that if $\alpha$ or $\beta$ equals $\pm\infty$, the interval should be interpreted as being unbounded in the associated direction(s).  It may be easily verified that  the set  $\cM$ is nonempty iff the following conditions hold:
\begin{equation}\label{nonempt}
\mu\in [\alpha,\beta]\;\;{\rm and}\;\;\sigma^2 \le (\beta-\mu)(\mu-\alpha),
\end{equation}
which will be assumed throughout. (We assume here that $0 \times \infty=0$, so that if, e.g.,  $\mu=\alpha$ and $\beta=+\infty$, then the right hand side of (\ref{nonempt}) is 0.)

Furthermore, one can also identify conditions under which $\cM$ is a singleton.

\begin{observation}\label{singletonian}
If $-\infty<\alpha<\beta <+ \infty$,  $\mu \in [\alpha,\beta]$, and $\sigma^2 = (\beta-\mu)(\mu-\alpha)$, then $\cM$ consists of the single probability measure which assigns to the point  $\alpha$ probability $p=\frac{\beta - \mu}{\beta - \alpha}$, and to the point  $\beta$ probability $1-p= \frac{\mu-\alpha}{\beta - \alpha}$.
\end{observation}

We now rephrase $\psi(x)$ as the optimal value of a certain optimization problem.  For use in later proofs, we define the following more general maximization problem, in terms of a general integrable objective function $\zeta$:
\begin{equation}\label{dual-1}
\begin{array}{clll}
 \sup\limits_{Q \in \cP(\ii)}&\int \zeta(\tau) dQ(\tau)\\
 {\rm s.t.}& \int \tau dQ(\tau)  =  \mu,\;
  \int \tau^2 dQ(\tau)  =  \mu^2 + \sigma^2.
\end{array}
\end{equation}
Our definitions imply that for all $x \in \bbr$, $\psi(x)$ equals the optimal value of Problem (\ref{dual-1}) for the special case that $\zeta(\tau) = \Psi(x,\tau)$.  Problem (\ref{dual-1}) is a classical  problem of moments \cite[see, e.g.,][]{lan:87}. From the Richter-Rogosinski Theorem \cite[e.g.,][Proposition 6.40]{SDR} or results in \cite{Bertsimas05},  we have the following.

\begin{observation}\label{rem-rr}
If Problem {\rm (\ref{dual-1})} possesses  an optimal solution, then it has an optimal solution with support of at most three points.
\end{observation}

\subsubsection{Review of relevant duality theory.}
As several of our later proofs will be based on duality theory, we now briefly review duality for Problem (\ref{dual-1}).  The  dual of Problem  (\ref{dual-1}) can be constructed as follows \cite[cf.][]{isii}.  Consider the Lagrangian
\[
L(Q,\lambda) :=  \int \Big[\zeta(\tau) - \sum_{i=0}^2 \lambda_i  \tau^i \Big]d Q(\tau)
+ \lambda_0 + \lambda_1 \mu + \lambda_2 (\mu^2 + \sigma^2).
\]
By maximizing $L(Q,\lambda)$ with respect to $Q \in \cP(\ii)$, and then minimizing with respect to $\lambda$, we obtain the following Lagrangian dual for Problem (\ref{dual-1}):
\begin{equation}\label{dual-2}
\begin{array}{clll}
\inf\limits_{\lambda\in \bbr^{3}}& \lambda_0+ \lambda_1 \mu + \lambda_2 (\mu^2 + \sigma^2)\\
{\rm s.t.}& \lambda_0 + \lambda_1 \tau + \lambda_2 \tau^2 \ge \zeta(\tau),\;\;\tau\in\ii.
\end{array}
\end{equation}
We denote by $\val(P)$ and $\val(D)$ the respective optimal values of the primal Problem (\ref{dual-1}) and its dual Problem (\ref{dual-2}).  By convention, if Problem (\ref{dual-1}) is infeasible, we set $\val(P) = -\infty$, and if Problem (\ref{dual-2}) is infeasible, we set $\val(D) =+ \infty$.  We denote by $\sol_P(x)$ the set of optimal solutions of the primal problem, and by $\sol_D(x)$ the set of optimal solutions of the dual problem, and note that these sets may be empty, even when both programs are feasible, e.g. if the respective optimal values are approached but not attained.

Note that  $\val(D) \ge \val(P)$.  We now give  sufficient conditions for there to be  no duality gap, i.e. $\val(P) = \val(D)$, as well as conditions for Problems (\ref{dual-1}) and (\ref{dual-2}) to have optimal solutions.  By specifying known general results for duality of such programs, e.g.,   \cite[Theorem 5.97]{BS}, to the considered setting, we have the following.

\begin{proposition}\label{dualprop1}
\label{pr-dual}
If $\bar{Q}$ is a probability measure which is feasible for the primal Problem {\rm (\ref{dual-1})}, $\bar{\lambda}=(\bar{\lambda}_0,\bar{\lambda}_1,\bar{\lambda}_2)$ is a vector which is feasible for the dual Problem {\rm (\ref{dual-2})}, and
\begin{equation}\label{stricter2}
\supp(\bar{Q}) \subseteq \big\lbrace \tau \in \ii : \zeta(\tau) = \bar{\lambda}_0 +  \bar{\lambda}_1 \tau + \bar{\lambda}_2 \tau^2  \big\rbrace,
\end{equation}
then $\bar{Q}$ is an optimal primal solution, $\bar{\lambda}$ is an optimal dual solution, and $\val(P) = \val(D)$.
Conversely, if $\val(P) = \val(D)$, and $\bar{Q}$ and $\bar{\lambda}$ are optimal solutions of the respective primal and dual problems, then condition {\rm (\ref{stricter2})} holds.
\end{proposition}

\subsubsection{Scarf's solution.}
We note that the distributionally robust single-stage news vendor problem considered by \cite{Scarf} is exactly Problem (\ref{nvc2}), when  $\ii = \bbr_+$.  As it will be useful for later proofs, we briefly review Scarf's explicit solution.  We actually state a slight generalization of the results of Scarf, and for completeness we include a proof in the technical appendix (Section\ \ref{appsec}). Define $\cf(z) := \big( (z - \mu)^2 + \sigma^2 \big)^{\frac{1}{2}}$ for all $z\in\bbr$.

\begin{theorem}\label{Scarfold}
Suppose that $b > c$,  $c + h > 0$, $\mu > 0$, $\sigma > 0$, and $\ii = \bbr_+$.  Let $\kappa := \frac{b - h - 2c}{b + h}$.  Then for each $x\in \bbr$,
\begin{equation}
\begin{aligned}
    \psi(x) = \begin{cases}  c\mu+\frac{b+h}{2}\big((x-\mu)^2+\sigma^2\big)^{\frac{1}{2}}
   -\frac{b-h-2c}{2}(x-\mu),& \text{\rm if}\ x\geq \frac{\mu^2+\sigma^2}{2\mu},\\
      \frac{(h+c)\sigma^2-(b-c)\mu^2}{\mu^2+\sigma^2}x+b\mu, & \text{\rm if}\  x \in[0, \frac{\mu^2+\sigma^2}{2\mu}),\\
      b \mu - (b - c) x, & \text{\rm otherwise}.
      \end{cases}
\end{aligned}
\end{equation}

As a consequence, a complete solution to the problem $\inf_{x \in \bbr} \psi(x)$ is as follows.
\begin{description}
\item [{\rm (i)}] If $\frac{\sigma^2}{\mu^2}>\frac{b-c}{h+c}$, then the unique optimal solution is $x = 0$, and the optimal value is $\mu b$.
\item [{\rm (ii)}] If
$\frac{\sigma^2}{\mu^2}<\frac{b-c}{h+c}$, then the unique optimal solution is $x = \mu + \kappa \sigma (1 - \kappa^2)^{-\frac{1}{2}}$, and the optimal value is $c \mu + \big( (h + c)(b - c) \big)^{\frac{1}{2}} \sigma$.
\item [{\rm (iii)}]
 If $\frac{\sigma^2}{\mu^2}=\frac{b-c}{h+c}$, then all $x \in [0,  \mu + \kappa \sigma (1 - \kappa^2)^{-\frac{1}{2}}]$ are optimal, and the optimal value is $\mu b$.
\end{description}
Furthermore,  in all cases $\argmax_{Q\in \cM}\bbe_Q [\Psi(x,D)]$ is nonempty  for every $x \in \bbr$.  Also, the optimal solution set and value of the problem $\inf_{x \in \bbr} \psi(x)$ is identical to that of Problem {\rm (\ref{nvc2})}, i.e. optimizing over $x \in \bbr$, as opposed to $x \in \bbr_+$, makes no difference.
\end{theorem}

For use in later proofs (e.g., sufficient conditions for strong time consistency in Theorem 4), it will also be useful to demonstrate a particular variant of Theorem\ \ref{Scarfold}.  Suppose that in Problem (\ref{nvc2}), we are not forced to select a deterministic constant $x$, but can instead select any distribution $D_1$ for $x$.  Specifically, let us consider the following minimax problem:
\begin{equation}\label{ddminmax}
\inf_{Q_1 \in \cP(\ii)} \phi(Q_1),
\end{equation}
where
$$\phi(Q_1) := \sup_{ Q_2 \in \cM} \bbe_{Q_1\times Q_2}\big[\Psi(D_1,D_2)\big],$$
and the notation $\bbe_{Q_1\times Q_2}$ indicates that for any choices for the marginal distributions $Q_1,Q_2$ of $D_1$ and $D_2$, the expectation is taken with respect to the associated product measure, under which $D_1$ and $D_2$ are independent.
In this case, we have the following result, whose proof we defer to the technical appendix (Section\ \ref{appsec}).

\begin{proposition}\label{Scarfold2}
Suppose that $b > c$,  $c + h > 0$, $\mu > 0$, $\sigma > 0$, $\frac{\sigma^2}{\mu^2} > \frac{b-c}{h+c}$, and
$\ii = \bbr$.  Then Problem {\rm (\ref{ddminmax})} has the unique optimal solution $\bar{Q}_1=\delta_0$.
\end{proposition}

We also note that $\psi$ inherits the property of convexity from $\Psi$.
\begin{observation}\label{convex1}
$\Psi(\cdot, d)$ is a convex function for every fixed $d\in \ii$, $\psi$ is a convex function on $\bbr$, and Problem {\rm (\ref{nvc2})} is a convex program.
\end{observation}

\subsubsection{Generalization of \cite{Natarajan} to a class of convex, continuous, piecewise affine functions.}\label{extendscarfsec}

\ 
\cite{Scarf} gave an explicit solution for Problems (\ref{dual-1}) and (\ref{dual-2}) when $\ii = \bbr_+$, and $\zeta$ is a convex, continuous piecewise affine function with exactly two pieces, by explicitly constructing a feasible primal - dual solution pair satisfying the conditions of Proposition\ \ref{pr-dual} (details of this construction can be found in Section\ \ref{appsec}).  \cite{Natarajan} generalized Scarf's results to a class of convex, continuous, piecewise affine (CCPA) functions with three pieces.  
We now state the solution to a special case of the problems studied in \cite{Natarajan}, as we will need the solution to such problems for our later studies of time consistency. For completeness, we provide a proof in the technical appendix (Section\ \ref{appsec}).

\begin{theorem}\label{useme1}[\cite{Natarajan}]
Suppose that there exist $c_1 , c_2 > 0$ such that  $c_1 < c_2$, and $\zeta(d) = \max\{- d + c_1, 0, d - c_2\}$ for all $d \in \bbr$.  Let $\eta := \frac{1}{2} (c_1 + c_2)$, and recall that $\cf(z) := \big( (z - \mu)^2 + \sigma^2 \big)^{\frac{1}{2}}$.  Further suppose that $\sigma > 0$, $\ii = \bbr_+$,
$$\frac{1}{4}(2 \mu - 3 c_1 + c_2)(3 c_2 - c_1 - 2 \mu) \leq \sigma^2,$$
and $\eta - \cf(\eta)\geq 0$.  Then the unique optimal solution to the primal Problem {\rm (\ref{dual-1})} is the probability measure $Q$ having support at two points $h_1=\eta - \cf(\eta)$ and $h_2=\eta + \cf(\eta)$, with
\begin{equation}\label{eq-th2}
\begin{array}{ll}
Q(h_1) = \sigma^2 \left( \sigma^2 + \big( \eta - \cf(\eta) - \mu \big)^2 \right)^{-1},\;\;
Q(h_2) = 1 - Q(h_1).
\end{array}
\end{equation}
Also, the unique optimal solution to the dual Problem {\rm (\ref{dual-2})} is
\begin{equation}\label{eq-th2-dual}
\lambda_0 = \half\big( \eta^2 + (\eta-\mu)^2 + \sigma^2 \big) \cf^{-1}(\eta)+\frac{c_1-c_2}{2} ,\ \ \ \lambda_1 = - \eta \cf^{-1}(\eta) ,\ \ \ \lambda_2 = \half \cf^{-1}(\eta),
\end{equation}
where $\cf^{-1}(\eta)$ represents the reciprocal of $\cf(\eta)$.
\end{theorem}

\section{Multistage formulation}\label{multisec}

In this section, we study a multistage extension of the distributionally robust news vendor problem discussed in Section \ref{sec-disrob}.

\subsection{Classical formulation}
\label{sec-multclass}
We begin by giving a quick review of  the classical (i.e. non-robust) multistage news vendor  problem (also called inventory problem), and start by introducing some additional notations.
For a vector $z=(z_1,...,z_n) \in \bbr^n$
and $1\le i\le j \le n$, denote $z_{[i,j]}:=(z_i,...,z_j)$. In particular for $i=1$ we simply write $z_{[j]}$ for   the vector consisting of the first $j$ components of $z$, and set $z_{[0]} := \emptyset$.

We suppose that  there is a finite time horizon $T$, and a (random)  vector of demands $D  = (D_1,\ldots,D_T)$. By $d  = (d_1,\ldots,d_T)$ we usually denote a particular realization of the random vector $D$.  We assume that the components of random vector $D$  are {\em mutually independent}, and refer to this as the {\em stagewise independence} condition.
We now define the family of admissible policies $\Pi$ by introducing two families of functions, $\lbrace y_t,\; t = 1,\ldots,T \rbrace$ and $\lbrace  x_t,\; t = 1,\ldots,T \rbrace$.  Conceptually, $y_t$ will correspond to the inventory level at the start of stage $t$, and $x_t$ will correspond to the inventory level after having ordered in stage $t$, but before the demand in that stage is realized.

We will consider policies which are nonanticipative, i.e. decisions do not depend on realizations of future demand.  We assume that $y_1$, the initial inventory level, is a given constant.  We also require that one can only order a nonnegative amount of inventory at each stage.
Thus the set of admissible policies $\Pi$ should consist of those vectors of (measurable) functions $\pi = \{x_t(d_{[t - 1]}), t = 1,\ldots,T\},$ such that  $x_t: \bbr_+^{t-1} \rightarrow \bbr$ satisfies $x_t(d_{[t-1]}) \ge y_t$, for all $d_{[t-1]} \in \bbr_+^{t - 1}$ and $t=1,...,T$, where
\begin{equation}\label{dyneqy}
y_{t+1} = x_t(d_{[t-1]}) -  d_t,\;t=1,...,T-1.
\end{equation}

It follows that any given choice of $\pi \in \Pi$, along with the given $y_1$, completely determines the associated functions $y_1,\ldots,y_T$.  Sometimes we will explicitly express $x_t$ and $y_t$ as a function of the associated policy $\pi$ and demands $d_{[t]}$ with the notations $x^{\pi}_t(d_{[t - 1]})$ and $y^{\pi}_t(d_{[t - 1]})$; other times we will suppress this notation.  Combining the above, we can write the classical multistage news vendor problem (inventory problem) as follows:

\begin{equation}
\label{invent-1classic}
\begin{aligned}
 & \inf_{\pi \in \Pi} \bbe
\left\{ \sum_{t=1}^T
 \rho^{t-1}\big[ c_t\big(x^{\pi}_t(D_{[t-1]}) - y^{\pi}_t(D_{[t-1]})\big)+ \Psi_t\big(x^{\pi}_t(D_{[t-1]}), D_t\big)\big] \right\}.
\end{aligned}
\end{equation}
Here $\rho \in (0,1]$ is a discount factor, $c_t, b_t, h_t$ are the ordering, backorder penalty and holding costs per unit in stage $t$, respectively, and
\begin{equation}\label{func-1}
\Psi_t(x_t,d_t):=   b_t[d_t-x_t]_+ + h_t[x_t-d_t]_+.
\end{equation}
Unless stated otherwise, we assume that $b_t>c_t>0$ and $h_t \geq 0$  for all $t=1,...,T$.

Problem (\ref{invent-1classic}) can be viewed as an optimal control problem in discrete time with state variables $y_t$, control variables $x_t$ and random parameters $D_t$.
It is well known that Problem (\ref{invent-1classic}) can be solved using DP equations, which can be written as
\begin{equation}\label{dyn-class-1}
 V_t(y_t)=\inf_{x_t \ge y_t}\left\{c_t(x_t - y_t)+\bbe\big[\Psi_t(x_t,D_t) + \rho V_{t+1}(x_t - D_t)\big]\right\},
\end{equation}
$t=1,...,T$, with $V_{T+1}(\cdot)\equiv 0$ (e.g., \cite{Zi}). Note that the value functions $V_t(\cdot)$ are convex, and do not depend on the demand data because of the stagewise independence assumption.
These equations naturally define a set of policies through the relation $x_t(y_t) \in \cX_t(y_t)$, where $\cX_t(y_t)$, $t=1,...,T,$  is the set of optimal solutions of the problem
\begin{equation}\label{dyn-class-2}
 \inf_{x_t \ge y_t }\big\{c_t(x_t - y_t)+\bbe[\Psi_t(x_t,D_t) + \rho V_{t+1}(x_t - D_t)]\big\},
\end{equation}
and the optimal value of Problem (\ref{invent-1classic}) is given by $V_1(y_1)$. Note that $x_t(y_t)$, $t = 1,\ldots,T$, are functions of $y_t$, i.e., it suffices to consider policies (measurable functions)  of the form $x_t=\pi_t(y_t)$; this fact is well known from optimal control theory (see, e.g., \cite{berts} for technical  details). Of course, the assumption of  stagewise  independence is essential for this conclusion.

Under the specified conditions,  the objective function of Problem (\ref{dyn-class-2}) tends to $+\infty$ as $x_t\to \pm\infty$. It thus follows from convexity that this objective function possesses a (possibly non-unique) unconstrained minimizer $x_t^*$ over $x\in \bbr$, and $\bar{x}_t:=\max\{y_t,x_t^*\}$ is an optimal solution of Problem (\ref{dyn-class-2}).  In particular, the so-called {\em base-stock  policy}  is optimal for the inventory Problem (\ref{invent-1classic}), where we note that such a result is classical in the inventory literature.

\begin{definition}\label{def-basest}
A policy $\pi \in \Pi$ is said to be a {\em base-stock} policy if there exist  constants $ x^*_t$, $t=1,\ldots,T,$ such that
 \begin{equation}\label{bstpol}
x^{\pi}_t  = \max\big\{y^{\pi}_t, x^*_t\big\},\; t = 1,\ldots,T,
 \end{equation}
\end{definition}

That is, Problem (\ref{invent-1classic}) can be solved using the DP formulation (\ref{dyn-class-1}) and associated policy (\ref{dyn-class-2}) in the following sense.

\begin{lemma}\label{classicsolve1}
The optimal value of Problem {\rm (\ref{invent-1classic})} equals $V_1(y_1)$.  Any policy $\pi$ such that $x^{\pi}_t(d_{[t - 1]}) \in\cX_t\big(y^{\pi}_t(d_{[t - 1]})\big)$ are for all   $t = 1,...,T$ and $d_{[t - 1]} \in \bbr_+^{t-1}$, is an optimal solution to Problem {\rm (\ref{invent-1classic})}.  Conversely, for any optimal policy $\pi$ for Problem {\rm (\ref{invent-1classic})}, and any $t \in \{1,...,T\}$, there exists a set $A \subseteq \bbr$ such that $\prob\big(y^{\pi}_t(D_{[t-1]}) \in A\big) = 1$, and $x^{\pi}_t(D_{[t-1]}) \in \cX_t\big(y^{\pi}_t(D_{[t-1]})\big)$ conditional on the event $\lbrace y^{\pi}_t(D_{[t-1]}) \in A \rbrace$.  Furthermore, it follows from the convexity of the relevant cost-to-go functions $V_t(y_t)$ that any set of base-stock constants $\lbrace x^*_t, t = 1,\ldots, T \rbrace$ such that
$x^*_t \in \cX_t(0)$ for all $t \in [1,T]$ will yield an optimal policy for Problem {\rm (\ref{invent-1classic})}.
\end{lemma}

As we shall see, such an equivalence does not necessarily hold for distributionally robust multistage inventory problems with moment constraints.

\subsection{Distributionally robust formulations}
\label{robustmulti}

Suppose now that the distribution of the demand process is not known, and we only have at our disposal information about the support and first and second order moments.  In this case, it is natural to use the framework previously developed for the single-stage problem (see  Section  \ref{singlestagesec}) to handle the distributional uncertainty at each stage.  However, in the multistage setting, there is a nontrivial question of how to model the associated uncertainty in the joint distribution of demand.

We will consider two formulations, one intuitively corresponding to the modeling choices of a policy maker who does not recompute her policy choices after each stage and one corresponding to a policy-maker who does. These two formulations are analogous to the two optimization models discussed in \cite{iye:05} and \cite{nil:05} in the framework of robust MDP, and can also be interpreted through the lens of (non)rectangularity of the associated families of priors (cf. \cite{ES,iye:05,nil:05}), as we will explore later in this section.  We refer to these formulations as {\em multistage-static} and distributionally robust DP, respectively.  Questions regarding the interplay between the sets of optimal policies of these two formulations are important from an implementability perspective, as a policy deemed optimal at time 0, but which does not remain optimal if the relevant decisions are re-examined at a later time, may not be implemented by the relevant stake-holders.  We note that such considerations were one of the original motivations for the study of time consistency in economics (cf. \cite{Strotz55}).  We further note that the particular definitions and formulations we introduce here are by no means the only way to define the relevant notions of time consistency, and again refer the reader to the survey by \cite{EJT}, and other recent papers in the optimization community (cf. \cite{iye:05,bod:06,carp,IPS,Homem.16,Shapiro17})    for alternative perspectives.

We suppose that we have been given a sequence of closed (possibly unbounded)  intervals $\ii_t = [\alpha_t,\beta_t] \subset  \bbr$, $t = 1,\ldots,T$, and sequences of the corresponding means $\lbrace \mu_t, t=1,\ldots,T \rbrace$, and variances  $\lbrace \sigma^2_t, t=1,\ldots,T \rbrace$.

\subsubsection{Multistage-static formulation.}\label{multistaticformsec}
We first consider the following  formulation, referred to as {\em multistage-static}, in which the policy maker does not recompute her policy choices after each stage.  Let us define
\begin{eqnarray}
 \label{probm-2}
\cM_t&:=&\left\{Q_t\in\cP(\ii_t):\bbe_{Q_t}[D_t]=
\mu_t,\;\bbe_{Q_t}[D_t^2]
 =\mu_t^2+\sigma_t^2 \right \},\;t=1,...,T;\label{last-2}\\
\label{probm-1}
 \cM&:=&\{Q=Q_1\times\cdots\times Q_T:Q_t\in \cM_t,\;t=1,...,T\}.
\end{eqnarray}
That is, the set $\cM$ consists of probability measures given by direct products of probability measures $Q_t\in \cM_t$. This can be viewed as an extension of the stage-wise independence condition, employed in Section \ref{sec-multclass}, to the considered distributionally robust case.
In  order for the sets $\cM_t$ to be nonempty we assume that (compare with (\ref{nonempt}))
\begin{equation}\label{nonem-2}
\mu_t\in [\alpha_t,\beta_t]\;\;{\rm and}\;\;\sigma_t^2 \le (\beta_t-\mu_t)(\mu_t-\alpha_t),\;t=1,...,T.
\end{equation}

According to (\ref{probm-1}), the associated minimax problem supposes that although the set of associated marginal distributions may be ``worst-case", the joint distribution will always be a product measure (i.e. the demand will be independent across stages).  The multistage-static formulation for the distributionally robust inventory problem can then be formulated as follows.
\begin{equation}\label{invent-1}
\begin{aligned}
 & \inf_{\pi \in \Pi}
 \sup_{Q\in \cM} \bbe_Q\big[
Z^\pi\big],
\end{aligned}
\end{equation}
where  $Z^\pi=Z^\pi (D_{[T]})$ is a function of
$D_{[T]}=(D_1,...,D_T)$ given by
\begin{equation}\label{inv-z}
Z^\pi (D_{[T]}):= \sum_{t=1}^T
 \rho^{t-1}\big[ c_t\big(x^{\pi}_t(D_{[t-1]}) - y^{\pi}_t(D_{[t-1]})\big)+ \Psi_t\big(x^{\pi}_t(D_{[t-1]}), D_t\big)\big],
\end{equation}
and $\Pi$ is the set of admissible policies defined previously in Section  \ref{sec-multclass}.
Of course, if the set $\cM=\{Q\}$ is a singleton, then formulation (\ref{invent-1}) coincides with formulation (\ref{invent-1classic}) taken with respect to the distribution $Q=Q_1\times\cdots\times Q_T$ of the demand vector $D_{[T]}$.

We note that the multistage-static formulation \eqref{invent-1} is closely related to optimization with risk measures. Indeed, the functional  
$\sup_{Q\in \cM} \bbe_Q\big[
Z\big]$ is a coherent risk measure (cf. \cite{SDR}).

Very little is known about the set of optimal policies for Problem\ (\ref{invent-1}), as this problem does not enjoy a DP formulation along the lines of (\ref{dyn-class-1}).

\subsubsection{Time consistency and distributionally robust DP equations.}

As informally referenced earlier, time inconsistency refers to the possibility that policy choices which seemed optimal from the perspective of time 0 no longer seem optimal if one re-performs one's minimax calculations at a later time.  Although first addressed within the economics community, the issue of time (in)consistency has recently started to receive attention in the stochastic and robust optimization communities (cf. \cite{Ri, bod:06, ADEHK,S-09,Ru,GH2011,carp,S-12,chen2013,IPS,Homem.16}), in which closely related concepts such as Pareto robust optimality (\cite{Iancu14}) have also been studied.  We note that related issues were addressed even in the seminal work of \cite{bell57} on DP, where it is asserted that:
  {\em ``An optimal policy has the property that whatever the initial state and initial decision are, the remaining decisions must constitute an optimal policy with regard to the state resulting from the first decision."}  The same principle has  been subsequently reformulated by several authors in a somewhat more precise form, e.g.,  in the recent work of \cite{carp}, where it is asserted that {\em
``The decision maker formulates
an optimization problem at time $t_0$ that yields a sequence of optimal decision
rules for $t_0$ and for the following time steps $t_1,...,t_N=T$. Then, at the next time
step $t_1$, he formulates a new problem starting at $t_1$ that yields a new sequence of
optimal decision rules from time steps $t_1$ to $T$. Suppose the process continues until
time $T$ is reached. The sequence of optimization problems is said to be dynamically
consistent if the optimal strategies obtained when solving the original problem at
time $t_0$ remain optimal for all subsequent problems."}  A nearly identical concept, which the authors refer to as the inherited-optimality-property (IOP), is also formalized in \cite{Homem.16}.

To motivate the particular definition of time-consistency we will use here, which is very similar to that used in e.g. \cite{Homem.16}, let us reason as follows.  Suppose we wished to know whether a policy $\pi$ which was optimal for the multi-stage static formulation had the property that, should one re-perform one's minimax calculation in the final period, one would make the same ordering decision.  As she cannot change past decisions, the only policy decision she still has to make is the determination of the function $x_T$.  However, she now has knowledge of $D_{[T-1]}$ and $y_T$, which she can incorporate into her minimax computations.  We note that here we are faced with the modeling question of how to reconcile the use of $D_{[T-1]}$ and $y_T$'s realized values in performing one's minimax computations with the previously assumed stagewise independence of demand.  A natural approach, consistent with the economics literature on time consistency, is to reason as follows.  As $D_{[T-1]}$ has already been realized, it is unreasonable to enforce independence of $D_T$ on this realization, as it is no longer undetermined.  Instead, the relevant minimax computation is carried out with this knowledge of the realization of $D_{[T-1]}$.  At that time, such a policy-maker is thus led to the optimization (at time T)
$$
\inf_{x_T \ge y_T}\left\{c_T(x_T - y_T) + \sup_{Q_T\in \cM_T}\bbe_{Q_T}[\Psi_T(x_T,D_T)] \right\},
$$
with 
$$
\cY_T(y_T) := \argmin_{x_T \ge y_T}\left\{c_T(x_T - y_T)+\sup_{Q_T\in \cM_T}\bbe_{Q_T}[\Psi_T(x_T , D_T)]\right\}
$$
the corresponding set of optimal policy choices.  Here we note that (for example) the inner maximization $\sup_{Q_T\in \cM_T}\bbe_{Q_T}[\Psi_T(x_T,D_T)]$ is implicitly a function of $D_{[T-1]}$, through the dependence on $x_T$.
Thus time-consistency of an optimal policy $\pi$ for the multi-stage static formulation should (at least as regards policy decisions in this final period) be equivalent to requiring that $x^{\pi}_T(D_{[T-1]}) \in \cY_T\big( y^{\pi}_T(D_{[T-1]}) \big).$  We note that there is a second subtlety here.  Indeed, as the exact distribution of $D_{[T]}$ is no longer known with certainty, the question of under which measures (for $D_{[T]}$) the requirement $x^{\pi}_t(D_{[t-1]}) \in \cY_T\big( y^{\pi}_T(D_{[T-1]}) \big)$ should hold w.p.1 must be resolved.  We note that when all distributions in ${\mathcal M}$ have the same support, such issues do not arise, while for the moment-based uncertainty sets we consider this distinction is important.  We also note that many past works do not take this subtlety into account in their definitions.  Here, we propose the natural and intuitive interpretation that one should require the inclusion hold w.p.1 for every measure in ${\mathcal M}$, as these are exactly those measures one believes possible.  

\textbf{Distributionally robust DP formulation.}  Before proceeding with our formal definition of time consistency, let us expand on the distributionally robust DP formulation, which we have defined only in the final period.  Carrying out the same logic inductively, we conclude that if a policy is to be deemed time-consistent when the policy-maker is (possibly) given the choice to recompute her minimax calculations in an arbitrary set of time periods, her choices should be consistent with the following distributionally robust DP equations.

\begin{equation}\label{dyn-mom-1}
 V_t(y_t)=\inf_{x_t \ge y_t}\left\{c_t(x_t - y_t)+\sup_{Q_t\in \cM_t}\bbe_{Q_t}[\Psi_t(x_t,D_t) + \rho V_{t+1}(x_t - D_t) ]\right\},
\end{equation}
$t=1,...,T$, with $V_{T+1}(\cdot)\equiv 0$.  The optimal value of DP formulation (\ref{dyn-mom-1}) is given by $V_1(y_1)$.
Dynamic equations (\ref{dyn-mom-1}) naturally define a set of policies of the form  $x_t=\pi_t(y_t)$,  $t = 1,\ldots,T$, with $x_t=\pi_t(y_t)$  being measurable selections $x_t\in \cY_t(y_t)$ from sets
 \begin{equation}\label{dyn-mompol}
 \cY_t(y_t) := \argmin_{x_t \ge y_t}\left\{c_t(x_t - y_t)+\sup_{Q_t\in \cM_t}\bbe_{Q_t}[\Psi_t(x_t , D_t) + \rho V_{t+1}(x_t - D_t)]\right\},\;t=1,...,T.
\end{equation}

We refer to \eqref{dyn-mom-1} as the \emph{distributionally robust DP formulation} and $V_1(y_1)$ as its optimal value.

We now observe that due to certain convexity properties, DP formulation (\ref{dyn-mom-1}) always possesses an optimal base-stock policy.  We note that such results are generally well-known to hold in this setting (cf. \cite{ACS}).  Recall Definition \ref{def-basest} of a base-stock  policy.

\begin{observation}\label{givesbound}
It follows from the convexity of the relevant cost-to-go functions $V_t(y_t)$ that Problem \eqref{dyn-mom-1} possesses an optimal base-stock policy.  Furthermore, any set of base-stock constants $\lbrace x^*_t, t = 1,\ldots, T \rbrace$ such that $x^*_t \in \cY_t(0)$ for all $t \in [1,T]$ yields an optimal policy.
Namely, for any such $\lbrace x^*_t, t = 1,\ldots, T \rbrace$, $\max\{y, x^*_t\} \in  \cY_t(y)$ for all $y \in \bbr$ and $t = 1,\ldots,T$.
\end{observation}

We note that the same conclusion could also have been drawn by rephrasing our formulation in the language of coherent risk measures, and applying known results for so-called nested risk measures (cf. \cite{rs06}, \cite[section 6.7.3]{SDR}), although we do not pursue such an analysis here.

We note that the question of whether or not there exists such an optimal base-stock policy for the multistage-static formulation is considerably more challenging, and will be central to our discussion of time consistency.

We close our discussion of the distributionally robust DP formulation with a final definition, formalizing our earlier discussion of for which measures one should require optimality (of decisions) under the distributionally robust DP formulation.
\begin{definition}[Robust-w.p.1-optimal]\label{def-robwp1}
Let us say that a policy $\pi \in \Pi$ is robust-w.p.1-optimal for the distributionally robust DP formulation if for all $Q \in {\mathcal M}$, w.p.1 $x^{\pi}_t(D_{[t-1]}) \in \cY_t\big( y^{\pi}_t(D_{[t-1]}) \big)$ for all $t \in [1,T]$.  
\end{definition}

\textbf{Formal definition of time consistency.}  We now formally define time consistency, in light of our earlier discussion.   We note that given the motivation behind time consistency, i.e. implementation of policies, a further subtlety must be considered.  Clearly, it is desirable for there to exist at least one policy which is optimal both initially, and if reconsidered at later times.  However, it is similarly undesirable for there to exist even one policy which could potentially be selected (i.e. optimal) initially, but deemed sub-optimal (i.e. non-implementable) at a later time.  Although such a notion is of course a stringent requirement (as noted informally in passing in \cite{Homem.16}), we believe that its conceptual importance none-the-less makes it worthy of further study.  This motivates the following definition(s) of time consistency, where we note that similar definitions were presented in \cite{GH2011} in a different context motivated by considerations in decision theory and artificial intelligence.  Then our definition of time consistency is as follows.

\begin{definition}[Time consistency]\label{def-consis}
If a policy $\pi \in \Pi$  is optimal for the multistage-static formulation  {\rm (\ref{invent-1})}, and robust-w.p.1-optimal for the distributionally robust DP formulation, we say that $\pi$ is {\rm time consistent}.
If there exists at least one optimal policy $\pi \in \Pi$ which is time consistent, we say that  Problem {\rm (\ref{invent-1})} is {\rm weakly time consistent}.
If every optimal policy of Problem {\rm  (\ref{invent-1})} is time consistent, we say that Problem {\rm (\ref{invent-1})} is {\rm strongly time consistent}.
\end{definition}

Of course the notion of strong time consistency makes sense only if Problem (\ref{invent-1}) possesses at least one optimal solution. Otherwise it is strongly time consistent simply because the set of optimal policies is empty.  

Our definition of time consistency can, in a certain sense, be viewed as an extension of the definition typically used in the theory of risk measures to an optimization context.  In Section\ \ref{diffvaluesex}, we show that it is possible for the multistage-static problem to have an optimal solution and to be strongly time consistent, but with a different optimal value than the distributionally robust DP formulation. That is, it is possible for the multistage-static problem to possess an optimal solution and to be strongly time consistent even when the rectangularity property does not hold.  This stands in contrast to the definition of consistency typically used in the theory of risk measures, i.e. the notion of dynamic consistency coming from \cite{ES} and based on a certain stability of preferences over time, which may result in a problem being deemed inconsistent based on the values that a given optimal policy takes under the different formulations, and even the values taken by suboptimal policies (cf. \cite{Ru,GH2011}).  In an optimization setting one may be primarily concerned only with the implementability of optimal policies, irregardless of their values and the values of suboptimal policies, and this is the approach we take here.  We note that such optimization-oriented formulations have been considered in several recent works, e.g. \cite{carp} and \cite{Homem.16}, and our definitions are largely consistent with those works.

Before exploring some of the subtle and interesting features of time (in)consistency for our model, we briefly review some previously known results for related models. Note that if the set $\cM$ is a singleton, then both the multistage-static and distributionally robust DP formulations collapse to the classical formulation. Hence both formulations have the same optimal value and strong time consistency follows. If one only has information about the support $\ii_t$, and hence takes $\cM_t$ to be the set of  all probability measures supported on the interval $\ii_t$, $t=1,...,T$, then both the multistage-static and distributionally robust DP formulations collapse to the so-called adjustable robust formulation (cf. \cite{BGGN}, \cite{dpadjust11}), which is purely deterministic. As a consequence, both formulations have the same optimal value and weak time consistency follows.  However, the recent work of \cite{Shapiro17} (itself inspired in part by an earlier version of this paper, \cite{Xin.13}), shows that even in that setting strong time consistency need not hold, and studies related phenomena in several particular problems with general moment constraints. However, we note that all inventory problems considered in \cite{Shapiro17} are weakly time consistent with both formulations having the same optimal value due to the rectangularity of the underlying set of measures.  \cite{Shapiro17} also shows that there are interesting connections between the notions of weak and strong time-consistency and the concept of ``strict monotonicity" for risk measures (e.g., \cite{Shapiro17a}), and we leave further investigations of this connection as an interesting direction for future research.  We also note that the existence of optimal time inconsistent policies was investigated earlier in several purely robust (i.e. deterministic) settings. In particular, \cite{BIP} demonstrated the optimality of so-called affine policies in certain settings, and \cite{Delage15} explicitly constructed optimal time inconsistent policies in an inventory control setting.

\textbf{Connection to rectangularity.}  To contextualize our definitions within the broader literature, we here briefly review the relevant notion of rectangularity.  Our definitions will closely follow those given in \cite{sha-14}, although we note that many closely related definitions have appeared previously throughout the literature (see e.g. \cite{iye:05}).  Consider the cost  $Z^\pi=Z^\pi (D_{[T]})$ of  a policy $\pi$, defined in (\ref{inv-z}).  Let $\widehat{\cM}$ be a set of probability distributions for the demand vector $D_{[T]}$, and let $Q\in \widehat{\cM}$. At the moment we do not assume that $Q$ is of the product form $Q=Q_1\times \cdots\times Q_T$, we will discuss this later. We can write

\begin{equation}
\label{eq-exp}
  \bbe_Q[Z^\pi]=
\bbe_{Q}\left[\bbe_{Q|D_1}\Big [\,\cdots\,  \bbe_{Q|D_{[T-2]}}\big [\bbe_{Q|D_{[T-1]}}[Z^\pi]\big]\Big]\right],
\end{equation}
where
$\bbe_{Q|D_{[t]}}[Z^\pi]$ is the conditional expectation, given $D_{[t]}$, with respect to the distribution $Q$ of $D_{[T]}$. Of course, this conditional expectation  is  a function of $D_{[t]}$.  Consequently,
\begin{equation}\label{expd2}
\sup_{Q\in \widehat{\cM}}\bbe_Q[Z^\pi]\le \sup_{Q \in \widehat{\cM}} \bbe_{Q }\left[\sup_{Q \in \widehat{\cM}}\bbe_{Q |D_1}\Big [\,\cdots\,  \sup_{Q \in \widehat{\cM}}\bbe_{Q |D_{[T-1]}}[Z^\pi]\Big]\right].
\end{equation}
The right hand side of (\ref{expd2}) leads to the nested formulation
\begin{equation}\label{invent-2b}
\inf_{\pi\in \Pi}\left\{ \sup_{Q \in \widehat{\cM}} \bbe_{Q}\left[\sup_{Q \in \widehat{\cM}}\bbe_{Q |D_1}\Big [\,\cdots\,  \sup_{Q \in \widehat{\cM}}\bbe_{Q |D_{[T-1]}}[Z^\pi]\Big]\right]\right\}.
\end{equation}
In particular, if the set $\widehat{\cM}$ is defined in the form (\ref{probm-1}), i.e., consists of products of probability measures (with the $t$-th measure drawn from $\cM_t$), then formulation  (\ref{invent-2b}) simplifies to
\begin{equation}\label{invent-d}
\inf_{\pi\in \Pi}\left\{ \sup_{Q_1\in \cM_1} \bbe_{Q_1}\left[\sup_{Q_2\in \cM_2}\bbe_{Q_2|D_1}\Big [\,\cdots\,  \sup_{Q_T\in \cM_T}\bbe_{Q_T|D_{[T-1]}}
[Z^\pi]\Big]\right]\right\}.
\end{equation}
It follows from \eqref{expd2} that the optimal value of (\ref{invent-2b}) is greater than or equal to the optimal value of the multistage-static Problem  (\ref{invent-1}). 
Moreover, the optimal value of (\ref{invent-2b}) can be strictly greater than the optimal value of (\ref{invent-1}). Let us demonstrate this through the following simple example.
\begin{example}\label{ex-1}
Let  $T=2$ and the set $\widehat{\cM}$ be of the product form (\ref{probm-1}). Suppose further  that $\I_1=[0,1]$, $\mu_1=1/2$ and $\sigma_1^2=1/4$. Then 
$\cM_1=\{Q_1\}$ is a singleton with $Q_1=p_1 \delta_0+p_2 \delta_1$,   $p_1=p_2=1/2$, i.e., with probability $1/2$ the demand $D_1$ can be either zero or one. Let us fix some policy $\pi \in \Pi$ and let  $Z^{\pi} = Z^{\pi}(D_1,D_2)$ be the corresponding objective function.  Then the associated cost under formulation (\ref{invent-1}) equals
\begin{equation}\label{stform-a}
   \sup_{Q_2\in \cM_2} \bbe_{Q_1\times Q_2}[Z^{\pi}(D_1,D_2)] =   \sup_{Q_2\in \cM_2} \Big(p_1\bbe_{Q_2}[Z^{\pi}(0,D_2)] +
p_2\bbe_{Q_2}[Z^{\pi}(1,D_2)] \Big),
\end{equation}
while the associated cost under formulation (\ref{invent-2b}) equals
\begin{equation}\label{stform-2a}
  \bbe_{Q_1}\left[\sup\limits_{Q_2\in \cM_2}\bbe_{Q_2|D_1}[Z^{\pi}(D_1,D_2)]\right]=
p_1   \sup\limits_{Q_2\in \cM_2}\bbe_{Q_2}[Z^{\pi}(0,D_2)]+
  p_2  \sup\limits_{Q_2\in \cM_2}\bbe_{Q_2}[Z^{\pi}(1,D_2)] .
\end{equation}
Note that the worst case distribution $Q_2\in \cM_2$ in (\ref{stform-a}) has to be the same for all possible realizations of the demand $D_1$.  In contrast, the worst case distribution $Q_2\in \cM_2$ in (\ref{stform-2a}) is allowed to depend on realized $D_1$. Hence the right hand side of (\ref{stform-2a}) can be strictly greater than the right hand side of (\ref{stform-a}).
\end{example}

In line with the definition given in \cite{sha-14}, we say that the set $\widehat{\cM}$ of probability measures is rectangular if such strict inequality does not occur for any r.v. $Z^{\pi}$, i.e. the two formulations are equivalent in their optimal values.  More formally, we make the following definition.

\begin{definition}
\label{def-rect}
Consistent with the definition given in \cite{sha-14}, we say that the set $\widehat{\cM}$ of probability measures is rectangular if for every measurable and non-negative function $f$, 
\begin{equation}\label{expd2bb}
\sup_{Q\in \widehat{\cM}}\bbe_Q[f(D_{[T]})] = \sup_{Q \in \widehat{\cM}} \bbe_{Q }\left[\sup_{Q \in \widehat{\cM}}\bbe_{Q |D_1}\Big [\,\cdots\,  \sup_{Q \in \widehat{\cM}}\bbe_{Q|D_{[T-1]}}[f(D_{[T]})]\Big]\right].
\end{equation}
\end{definition}
We note that under additional compactness assumptions on the support of $D_{[T]},$ \cite{sha-14} formally explores several related concepts and subtleties of this definition, e.g. proves that one can associate a rectangular set of probability measures to any given (possibly non-rectangular) set of probability measures, but as such a compactness condition does not hold in our setting we do not explore that further here.
\\\indent For a rectangular set $\widehat{\cM}$ the static formulation
\begin{equation}\label{expd4}
\inf_{\pi\in \Pi} \sup_{Q \in \widehat{\cM}} \bbe_{Q} \big[Z^\pi\big],
\end{equation}
is equivalent (in an appropriate sense, see e.g. \cite{sha-14,Shapiro17a}) to the formulation (\ref{invent-2b}).  Furthermore, the natural generalization of the distributionally robust DP equations (\ref{dyn-mom-1}) can be applied to (\ref{expd4}), with both formulations having a common optimal policy and the same optimal value.

We note that the concept of rectangularity has been central to the past literature on time consistency (cf. \cite{ES,GH2011,IPS}), especially as it relates to optimization (cf. \cite{iye:05,nil:05,WKR}). In several of these works, connections were made between tractability of the associated robust MDP and various notions of rectangularity (e.g. (s,a)-rectangularity, s-rectangularity). We refer the interested reader to \cite{WKR} and the references therein for details. Our definition of rectangularity is aimed directly at the decomposability property of the  static formulation ensuring its equivalence to the corresponding dynamic formulation (see \cite{sha-14} for details).

In general the set of product measures $\cM$ we consider in this work is not rectangular, as certified by the possible lack of weak time-consistency which we will soon demonstrate.  We note that a rectangular analogue of the set of measures $\cM$ defined in \eqref{probm-2}-\eqref{probm-1} would be the set of all joint distributions $Q$ for $D_{[T]}$ such that
\begin{eqnarray}
 \label{probm-2b}
D_t\in\cP(\ii_t)\ \ \ ,\ \ \ \bbe_{Q}[D_t | D_{[t-1]}]= \mu_t\ \ \ ,\ \ \ \bbe_{Q}[D_t^2| D_{[t-1]}] = \mu_t^2 + \sigma_t^2, \ \ \ t=1,\ldots,T.
\end{eqnarray}

Non-rectangular (and intractable) formulations for robust MDP are described in both \cite{iye:05} and \cite{nil:05}.  In \cite{iye:05}, it is referred to as the static formulation, while in \cite{nil:05}, it is referred to as the stationary formulation.  In both of these settings, these non-rectangular formulations essentially equate to requiring nature to select the same transition kernel every time a given state (and action, depending on the formulation) is encountered, as opposed to being able to select a different kernel every time a given state is visited in the robust MDP, and we refer the reader to \cite{iye:05}, \cite{nil:05}, and \cite{WKR} for details.  Although our multistage-static formulation could similarly be phrased in terms of a particular kind of dependency between the choices of nature in a robust MDP framework, and would be significantly different from either of the aforementioned non-rectangular formulations, we do not pursue such an investigation here, and leave the formalization of such connections as a direction for future research.

\section{Time consistency : sufficient conditions and (counter) examples}\label{consistencysection}\label{sec-tc1}

\subsection{Sufficient conditions for weak time consistency}\label{suffsecweak}
In this section, we provide simple sufficient conditions for the weak time consistency of Problem (\ref{invent-1}).  Our condition is essentially equivalent to monotonicity of the associated base-stock constants.  Intuitively, in this case the inventory manager can always order up to the optimal inventory level with which to enter the next time period, irregardless of previously realized demand.  Thus any potential for the adversary to take advantage of previously realized demand information in the distributionally robust DP formulation is ``masked" by the fact that the actual attained inventory level will always be this idealized level, under both formulations.  We note that several previous works have identified monotonicity of base-stock levels as a condition which causes various inventory problems to become tractable, in a variety of settings (cf. \cite{Veinott65,IV}, \cite{Jagannathan78}, \cite{Zi}).  In particular, \cite{Jagannathan78} studied a similar distributionally robust inventory model with moment constraints and identified  monotonicity of base-stock levels as a sufficient condition for a myopic base-stock policy to be optimal.  For completeness of the paper (as well as use in the later proofs), we state a variant of the results of \cite{Jagannathan78} in this section and include a proof in the appendix.

We begin by providing a different (but equivalent) formulation for Problem (\ref{invent-1}), in which all relevant instances of $y_t$ are rewritten in terms of the appropriate $x_t$ functions, as this will clarify the precise structure of the relevant cost-to-go functions.  As a notational convenience, let $c_{T+1} = 0$, in which case we define
\begin{equation}\label{inv-Psi}
\hat{\Psi}_t(x_t,d_t) := (c_t - \rho c_{t+1}) x_t + b_t[d_t - x_t]_+ +  h_t[x_t - d_t]_+,\;t=1,...,T.
\end{equation}
Let us define the problem
\begin{equation}\label{invent-redo1}
\begin{aligned}
 & \inf_{\pi \in \Pi}
 \sup_{Q\in \cM} \bbe_Q \left[ \sum_{t=1}^T \rho^{t-1} \hat{\Psi}_t\big(x_t(y_t),D_t\big) \right] - c_1 y_1 + \sum_{t=1}^{T-1} \rho^t c_{t+1} \mu_t.
\end{aligned}
\end{equation}
Then, using straightforward substitution we can make the following observation.
\begin{observation}\label{equalformulations}
Problem {\rm (\ref{invent-1})} and Problem {\rm (\ref{invent-redo1})} are equivalent, i.e. each policy $\pi \in \Pi$ has the same value under both formulations.
\end{observation}

We now derive a lower bound for any policy by allowing the policy maker to reselect her inventory at the start of each stage, at no cost.  As it turns out, this bound is ``attainable" when the set of base-stock levels is monotone increasing.  For $x \in \bbr$, let us define
\begin{equation}
\label{eq-defs}
\eta_t(x) := \sup\limits_{Q_t \in \cM_t} \bbe_{Q_t} [ \hat{\Psi}_t(x,D_t) ],\;\;
\Gamma^x_t :=  \argmax\limits_{Q_t \in \cM_t} \bbe_{Q_t} [ \hat{\Psi}_t(x,D_t) ],
\end{equation}
and let
\begin{equation}
\label{eq-defs-2}
\begin{array}{lll}
\hat{\eta}_t := \inf\limits_{x \in \bbr} \eta_t(x)
= \inf\limits_{x \in \bbr} \sup\limits_{Q_t \in \cM_t} \bbe_{Q_t} [ \hat{\Psi}_t(x,D_t) ],\\
\hat{\Gamma}_t := \argmin\limits_{x \in \bbr} \eta_t(x)=\argmin\limits_{x \in \bbr}
\sup\limits_{Q_t \in \cM_t} \bbe_{Q_t} [ \hat{\Psi}_t(x,D_t) ].
\end{array}
\end{equation}

For $j \geq 1$, and probability measures $Q_1,\ldots,Q_j$, let us define $\otimes_{t=1}^j Q_t :=Q_1 \times \ldots \times Q_j$, i.e. the associated product measure with the corresponding marginals.  Then we have the following.

\begin{lemma}\label{lboundme}
Suppose that the sets  $\Gamma^x_t,\hat{\Gamma}_t$ are non-empty for all $x \in \bbr$, $t =1,...,T$.  Let us fix any $\pi = (x_1,\ldots,x_T) \in \Pi$, and $i \geq 0$.  Then for any given $Q_1 \in \cM_1, \ldots, Q_i \in \cM_i$, there exist $Q_{i+1} \in \cM_{i+1},\ldots,Q_T \in \cM_T$ such that
\begin{equation}\label{inductme}
\bbe_{\otimes_{j=1}^T Q_j}\big[\hat{\Psi}_t\big(x_t(y_t),D_t\big)\big] \geq \hat{\eta}_t\ \textrm{for all}\ t \geq i + 1.
\end{equation}
Furthermore, the optimal value of Problem {\rm (\ref{invent-1})} is at least $\sum_{t=1}^T \rho^{t-1} \hat{\eta}_t  - c_1 y_1 + \sum_{t=1}^{T-1} \rho^t c_{t+1} \mu_t.$
\end{lemma}

We defer the proof to the technical appendix (Section\ \ref{appsec}). We now show that the bound of Lemma\ \ref{lboundme} is ``realizable" when the set of base-stock levels is monotone increasing, and that in this case the associated base-stock policy is optimal for both the multistage-static and distributionally robust DP formulations.  In particular, in this setting, the associated base-stock policy is time consistent, and thus the multistage-static problem is weakly time consistent. Again we defer the proof to the appendix Section\ \ref{appsec}.

\begin{theorem}\label{isbsest}
Suppose there exists a nondecreasing sequence $x^*_t,$ $t=1,...,T$, such that  $y_1 \leq x^*_1$,   and $x^*_t \in \hat{\Gamma}_t$, $t =1,...,T$, where $\hat{\Gamma}_t$ is defined in {\rm (\ref{eq-defs-2})}.  Also suppose $\ii_t \subset \bbr_+$ for all $t=1,...,T$.  Then the base-stock policy $\pi$ for which $x_t(y) = \max\{y,x^*_t\}$ for all $y \in \bbr$, is an optimal policy for the multistage-static formulation, and a robust-w.p.1-optimal policy for the distributionally robust DP formulations, and attains value
 $\sum_{t=1}^T \rho^{t-1} \hat{\eta}_t - c_1 y_1 + \sum_{t=1}^{T-1} \rho^t c_{t+1} \mu_t$ under both formulations.  Consequently, this base-stock policy is time consistent, and the multistage-static problem is weakly time consistent.
\end{theorem}

We note that Theorem \ref{isbsest} implies that if the parameters
$\mu_t,\sigma_t,c_t,b_t,h_t$ and $\ii_t$ are the same for all $t=1,...,T$, and hence the sets
$\cM_t$ are also the same for all $t$,
then the multistage-static problem is weakly time consistent, and the multistage-static and distributionally robust DP formulations have the same optimal value.

\subsection{Sufficient conditions for strong time consistency}\label{suffsecstrong}

In this section, we show that under additional assumptions, which ensure that the variance in each stage is sufficiently large, the multistage-static problem is strongly time consistent.  As we will see, in this case
there is a unique optimal base-stock policy, and in this policy all base-stock constants equal zero, the intuition being that when the variance is sufficiently large, it becomes undesirable to give nature any additional ``wiggle room".  We further note that such a base-stock policy has been widely adopted in practice and the resulting inventory system is a so-called Make-To-Order (MTO) or ``Pull" system. In such a system, no inventory is carried and the replenishment is based on actual demands instead of forecasts (cf. \cite{Williams84, Arreola-Risa98, Federgruen99, Rajagopalan02, Kaminsky09}).  We will later see in Section\ \ref{weaknotstrongex} that deviating slightly from this setting may lead to a lack of strong time consistency.  In particular, our results demonstrate that strong time consistency is a very fragile property. Our sufficient conditions are as follows.

\begin{theorem}\label{consistent}
Suppose that $b'_t := b_t - c_t + \rho c_{t+1} > 0$, $h'_t := h_t + c_t - \rho c_{t+1} > 0$, $\sigma_t,\mu_t > 0$, $\ii_t = \bbr_+, t = 1,\ldots,T$, $y_1 = 0$, and
\begin{equation}\label{sufcontc}
\frac{\sigma^2_t}{\mu^2_t} > \frac{b'_t}{h'_t},\;\;t=1,...,T.
\end{equation}
Then the set of optimal policies for the multistage-static problem is exactly the set of policies $$\Pi^0 := \big\lbrace \pi = (x_1,\ldots,x_T) \in \Pi : x_1(y_1) = 0, x_t(z) = 0\; for\;all\;z \leq 0\; and \;t \in [1,T] \big\rbrace,$$
and the multistage-static problem is strongly time consistent.
\end{theorem}
We defer the proof to the technical appendix (Section 6). 
\\\indent We note that under certain rectangularity-related assumptions, necessary and sufficient conditions for the existence of time-inconsistent optimal policies in an inventory setting with moment constraints was very recently provided in \cite{Shapiro17}.  However, those results are not applicable to the setting we consider, as the uncertainty sets we consider here are inherently non-rectangular, and thus (for example) our formulation allows for the non-existence of weak time-consistency, as well as the two formulations having different optimal values, and even the possibility that no policy of base-stock form is optimal for the multistage-static formulation (while the assumptions of \cite{Shapiro17} do not allow for such behavior).  Furthermore, \cite{Shapiro17} considers only 2-period problems, while the sufficient conditions provided in this work hold in the general multi-period setting.

\subsection{Further investigation of time (in)consistency}\label{moretimesec}
We now demonstrate that the question of time (in)consistency becomes quite delicate for inventory models with moment constraints, by considering a series of examples in which our model exhibits interesting (and sometimes counterintuitive) behavior.
In particular: (i) the problem can fail to be weakly time consistent, (ii) the problem can be weakly but not strongly time consistent, and (iii) the problem can be strongly time consistent even if every associated optimal policy takes different values under the multistage-static and distributionally robust DP formulations.  We also prove that, although the distributionally robust DP formulation always has an optimal policy of the base-stock form, there may be no such optimal policy for the multistage-static formulation.  We note that (i) and (ii) are subtle phenomena which the simpler models discussed in several previous works  (e.g. \cite{S-12}) cannot exhibit.  We also note that (iii) emphasizes  an interesting and surprising feature of our model and definitions: (strong) time consistency can hold even when the underlying family of measures from which nature can select is non-rectangular.  This stands in contrast to much of the related work on time consistency, where rectangularity is essentially taken as a pre-requisite for time consistency.  We also note that (iii) stands in contrast to some alternative, less policy-focused definitions of time consistency, e.g. those definitions appearing in the literature on risk measures (cf. \cite{ES}), under which time consistency could not hold if an optimal policy took different values under the two formulations.
We view our results as a step towards understanding the subtleties which can arise when taking a policy-centric view of time consistency in an operations management setting.  Throughout this section, we will let $\Pi_s^{opt}$ denote the set of all optimal policies for the corresponding multistage-static problem, and $\Pi_d^{opt}$ denote the set of all robust-w.p.1-optimal policies for the corresponding distributionally robust DP problem.

\subsubsection{Example: a multistage-static problem that is not weakly time consistent.}\label{notweakex}
In this section, we explicitly provide an example for which the multistage-static problem is not weakly time consistent.  Furthermore, for this example, the multistage-static and distributionally robust DP formulations have different optimal values.

Let us define $y_1 = 10$, $\rho = 1,$
$$\ii_1= [1,3] ,\ \ \ \mu_1 = 2,\ \ \ \sigma_1 = 1,\ \ \ \ c_1 = 0,\ \ \ b_1 = 2,\ \ \ h_1 = 2,$$
$$\ii_2 = \bbr_+,\ \ \ \mu_2 = 8,\ \ \ \sigma_2 = 2,\ \ \ c_2 = 0,\ \ \ b_2 = 1,\ \ \ h_2 = 1.$$
Let $\tilde{\Pi}_s$ denote the set of policies $\tilde{\pi} = (\tilde{x}_1,\tilde{x}_2)$ such that  $\tilde{x}_1(10) = 10$, $\tilde{x}_2(9) = 9$, $\tilde{x}_2(7) = 7$, and
$\tilde{\Pi}_d$ denote the set of policies $\tilde{\pi} = (\tilde{x}_1,\tilde{x}_2)$ such that  $\tilde{x}_1(10) = 10$, $\tilde{x}_2(9) = 9$, $\tilde{x}_2(7) = 8$.
Note that here (and in later statements) we have defined a set of policies by specifying a required behavior at only a few values, and allow the behavior at all other values to be arbitrary (subject to the overall policy belonging to $\Pi$).

\begin{theorem}\label{notweak1}
$\Pi_s^{opt} = \tilde{\Pi}_s$,
and the optimal value of the multistage-static problem is 18.  On the other hand, $\Pi_d^{opt} \subseteq \tilde{\Pi}_d$, and the optimal value of the distributionally robust DP problem is $17+\frac{\sqrt{5}}{2} > 18$.  Consequently, the multistage-static problem is not weakly time consistent, and the multistage-static and distributionally robust DP problems have different optimal values.
\end{theorem}
We defer the proof to the technical appendix (Section\ \ref{appsec}).

\subsubsection{Example: a multistage-static problem that is weakly time consistent, but \emph{not} strongly time consistent.}\label{weaknotstrongex}

In this section, we explicitly provide an example showing that it is possible for the multistage-static problem to be weakly time consistent, but not strongly time consistent.  In particular, there is a base-stock policy $\pi^*$, with associated base-stock constants $x^*_1,x^*_2$ satisfying the conditions of Theorem \ref{isbsest}, which is optimal for the multi-stage static formulation and robust-w.p.1-optimal for the distributionally robust DP formulation, yet the multistage-static problem has other non-trivial optimal policies which are not robust-w.p.1-optimal for the distributionally robust DP formulation.  The intuitive explanation is as follows.  In the multistage-static formulation, one can leverage the randomness in the realization of $D_1$ to construct a policy $\pi'$ such that  with positive probability $x^{\pi'}_2(y_2)$ is slightly below $x^*_2$, and with the remaining probability is slightly above $x^*_2$.  Since in the multistage-static formulation nature cannot observe the realized inventory in stage 2 before selecting a worst-case distribution, it turns out that such a policy incurs the same cost as $\pi'$ under the multistage-static formulation.  However, under the distributionally robust DP formulation, such a perturbation leads to sub-optimality.  We note that such a lack of strong time consistency can also be interpreted as resulting from the fact that optimality of a policy for the static formulation does not require optimality for every possible measure which nature can select, analogous to the ideas explored (in the robust optimization setting) in \cite{Iancu14}.  We note that in this example, even though the multistage-static problem is not strongly time consistent, both formulations have the same optimal value, as dictated by Theorem \ref{isbsest}.

Let us define $y_1 = 0$, $\rho = 1,$
$$\ii_1= [1,3],\ \ \ \mu_1 = 2,\ \ \ \sigma_1 = 1,\ \ \ \ c_1 = 0,\ \ \ b_1 = 1 ,\ \ \ h_1 = 1,$$
$$\ii_2 = \bbr_+ ,\ \ \ \mu_2 = 10 ,\ \ \ \sigma_2 = 1,\ \ \ c_2 = 0 ,\ \ \ b_2 = 1,\ \ \ h_2 = 1.$$
Then we prove the following, whose proof we defer to the technical appendix (Section\ \ref{appsec}).
\begin{theorem}\label{weaknotstrong}
The multistage-static problem is weakly time consistent, but not strongly time consistent.
\end{theorem}

\subsubsection{Example: a multistage-static problem that is strongly time consistent, but the two formulations have a \emph{different} optimal value.}\label{diffvaluesex}
In this section, we explicitly provide an example showing that it is possible for the multistage-static problem to be strongly time consistent, yet for the two formulations to have different optimal values.  We note that, although it is expected that there will be settings where the two formulations have different optimal values, it is somewhat surprising that this is possible even when the two formulations have the same set of optimal policies.  As discussed previously, we note that this possibility stands in contrast to several related works which consider alternative, less policy-focused definitions of time consistency, e.g. those definitions appearing in the literature on risk measures.

Let us define $y_1 = 0$, $\rho = 1,$
$$\ii_1= [1,3] ,\ \ \ \mu_1 = 2,\ \ \ \sigma_1 = 1,\ \ \ \ c_1 = 0,\ \ \ b_1 = 0,\ \ \ h_1 = 0,$$
$$\ii_2 = \bbr_+,\ \ \ \mu_2 = 100,\ \ \ \sigma_2 = 5,\ \ \ c_2 = 2,\ \ \ b_2 = 1,\ \ \ h_2 = 1.$$
Let $\tilde{\Pi}$ denote the set of policies $\tilde{\pi} = (\tilde{x}_1,\tilde{x}_2)$ such that  $\tilde{x}_1(0) = 102$, $\tilde{x}_2(101) = 101$, $\tilde{x}_2(99) = 99$.  Then we prove the following,  whose proof we defer to the technical appendix (Section\ \ref{appsec}).

\begin{theorem}\label{strongdiff}
$\Pi_s^{opt} = \tilde{\Pi}$,
and the multistage-static problem is strongly time consistent.  However, the optimal value of the multistage-static problem equals 5, while the optimal value of the distributionally robust DP problem equals $\sqrt{26} > 5$.
\end{theorem}

\subsubsection{Example: a multistage-static problem that has no optimal policy of base-stock form.}\label{nobaseex}
In this section, we explicitly provide an example showing that it is possible for the multistage-static problem to have no optimal base-stock policy, where we note that in all our previous examples the associated multistage-static problem did indeed have an optimal base-stock policy (possibly different from that of the associated distributionally robust DP problem).  Note that this stands in contrast to the distributionally robust DP formulation, which always has an optimal base-stock policy by Observation\ \ref{givesbound}.  It remains an interesting open question to develop a deeper understanding of the set of optimal policies for the multistage-static problem, where we again note that some preliminary investigations of such distributionally robust problems with independence constraints can be found in \cite{GL2013}.  Both the results of \cite{GL2013}, and our own result, indicate that the structure of optimal policies for the multistage-static problem may be very complicated.
\\\indent To prove the desired result, it will be useful to consider a family of problems parameterized by a parameter $\epsilon$.  In particular, let $\epsilon \in \big(0, \frac{1}{2}(\sqrt{6} - 2) \big)$ be any sufficiently small strictly positive number.  It may be easily verified that for any such $\epsilon$, one has $\epsilon \in (0,\frac{1}{4})$, and
\begin{equation}\label{tc-epsilon}
\frac{1}{2} - 2\epsilon - \epsilon^2 > 0.
\end{equation}
Let us define $y_1 = 10 - \epsilon$, $\rho = 1,$
$$\ii_1= [1 - \epsilon,\ 3+\epsilon] ,\ \ \ \mu_1 = 2,\ \ \ \sigma_1 = 1,\ \ \ \ c_1 =
0,\ \ \ b_1 = 2,\ \ \ h_1 = 2,$$
$$\ii_2 = \bbr_+,\ \ \ \mu_2 = 8,\ \ \ \sigma_2 = 3,\ \ \ c_2 = 0,\ \ \ b_2 = 1,\ \
\ h_2 = 1.$$
Then we prove the following, whose proof we defer to the technical appendix (Section\ \ref{appsec}).
\begin{theorem}\label{nobase-stock}
Suppose $\epsilon$ satisfies \eqref{tc-epsilon}. Then any admissible policy
$\tilde{\pi}=(\tilde{x}_1,\tilde{x}_2)\in \Pi$ satisfying $\tilde{x}_1(y_1)=y_1$,
$\tilde{x}_2(D_1) = y_1 - D_1 + \epsilon$ belongs to $\Pi_s^{opt}$, and the corresponding optimal value equals
$19 - 2\epsilon$. Moreover, no base-stock policy belongs to $\Pi_s^{opt}$.
\end{theorem}

\section{Conclusion}\label{concsec}
In this paper, we analyzed the notion of time consistency in the context of managing an inventory under distributional uncertainty.  In particular, we studied the associated multistage distributionally robust optimization problem, when only the mean, variance  and distribution support are known for the demand at each stage.  Our contributions were three-fold.  First, we gave a refined policy-centric definition for time consistency in this setting, and put our definition in the broad context of prior work on time consistency and rectangularity.  More precisely, we defined two natural formulations for the relevant optimization problem.  In the multistage-static formulation, the policy-maker cannot recompute her policy after observing realized demand.  In the distributionally robust DP formulation, she is allowed to reperform her minimax computations at each stage.  If there exists a policy which is optimal for both formulations (w.p.1 under every joint distribution for demand belonging to the uncertainty set), we say that the policy is \emph{time consistent}, and the problem is \emph{weakly time consistent}.  If every optimal policy for the multistage-static formulation is time consistent, we say that the problem is \emph{strongly time consistent}.  
\\\indent Second, we gave sufficient conditions for weak and strong time consistency.  Intuitively, our sufficient condition for weak time consistency coincides with the existence of an optimal base-stock policy in which the base-stock constants are monotone increasing.  Our sufficient condition for strong time consistency can be interpreted in two ways.  On the one hand, strong time consistency holds if the unique optimal base-stock policy for the distributionally robust DP formulation is to order-up to 0 at each stage, i.e., the well-known Make-To-Order policy.  Alternatively, we saw that this condition also has an interpretation in terms of requiring that the demand variances are sufficiently large relative to their respective means.  Third, we gave a series of examples of two-stage problems exhibiting interesting and counterintuitive time (in)consistency properties, showing that the question of time consistency can be quite subtle in this setting.  In particular: (i) the problem can fail to be weakly time consistent, (ii) the problem can be weakly but not strongly time consistent, and (iii) the problem can be strongly time consistent even if every associated optimal policy takes different values under the multistage-static and distributionally robust DP formulations.  We also proved that, although the distributionally robust DP formulation always has an optimal policy of base-stock form, there may be no such optimal policy for the multistage-static formulation.  This stands in contrast to the analogous setting, analyzed in \cite{S-12}, in which only the mean and support of the demand distribution is known at each stage, for which it is known that such phenomena cannot occur (as the problem is always weakly time consistent).  
\\\indent We departed from much of the past literature by demonstrating both negative \emph{and positive} results regarding time consistency when the underlying family of distributions from which nature can select is non-rectangular, a setting in which most of the literature focuses on demonstrating hardness of the underlying optimization problems and other negative results.  Furthermore, our example demonstrating that it is possible for the multistage-static problem to be strongly time consistent, but with a different optimal value than the distributionally robust DP formulation, stands in contrast to the definition of time consistency typically used in the theory of risk measures, i.e. the notion of dynamic consistency coming from \cite{ES}, under which a problem may be deemed time inconsistent based on the values that a given optimal policy takes under the different formulations, and even the values taken by suboptimal policies.  Indeed, our definitions are motivated by the fact that in an optimization setting, one may be primarily concerned only with the implementability of optimal policies, irregardless of their values and the values of suboptimal policies, building on the more optimization-oriented definitions provided in \cite{carp} and the recent work \cite{Homem.16}.
\\\indent Our work leaves many interesting directions for future research.  The general question of time consistency remains poorly understood.  Furthermore, our work has shown that this question can be quite subtle.  For the particular model we consider here, it would be interesting to develop a better understanding of precisely when time consistency holds.  It is also an intriguing question to understand how much our two formulations can differ in optimal value and policy, even when time inconsistency occurs, along the lines of \cite{Huang11},  \cite{Asamov15}, and \cite{IPS}.  On a related note, it is largely open to develop a broader understanding of the optimal solution to the multistage-static problem, or even approximately optimal solutions, as well as related algorithms, where we note that preliminary investigations along these lines were recently carried out in \cite{GL2013}.  Of course, it is also an open challenge to understand the question of time consistency more broadly, how precisely the various definitions of time consistency presented throughout the literature relate to one-another, and more generally to understand the relationship between different ways to model multistage optimization under uncertainty.

\bibliographystyle{nonumber}

\begin{thebibliography}{}




\bibitem[{Ahmed, Cakmak and Shapiro(2007)}]{ACS}
Ahmed, S., U. Cakmak, A. Shapiro. 2007.
Coherent risk measures in inventory problems. {\it European Journal of Operational Research}. {\bf 182} 226--238.


\bibitem[{Arreola-Risa and DeCroix(1998)}]{Arreola-Risa98}
Arreola-Risa, A., G.A. DeCroix. 1998.
Make-to-order versus make-to-stock in a production-inventory system with general production times. {\it IIE Transactions}. \textbf{30}(8) 705-713.


\bibitem[{Artzner et al.(2007)}]{ADEHK}
Artzner, P., F. Delbaen, J.M. Eber, D. Heath, H. Ku. 2007.
Coherent multiperiod risk adjusted values and Bellmans principle. {\it Annals of Operations Research}. {\bf 152} 5--22.


\bibitem[{Asamov and Ruszczy\'{n}ski(2015)}]{Asamov15}
Asamov, T., A. Ruszczy\'{n}ski. 2015.
Time-consistent approximations of risk-averse multistage stochastic
optimization problems. {\it Mathematical Programming}.  {\bf 153}(2) 459-493.




\bibitem[{Bellman(1957)}]{bell57}
Bellman, R.E. 1957.
{\em Dynamic Programming}. Princeton University Press, Princeton, NJ.




\bibitem[{Ben-Tal et al.(2005)}]{BGNV}
Ben-Tal, A., B. Golany, A. Nemirovski, J.P. Vial. 2005.
Retailer-supplier flexible commitments contracts: a robust optimization approach. {\it Manufacturing \& Service Operations Management} {\bf 7} 248--271.


\bibitem[{Ben-Tal et al.(2004)}]{BGGN}
Ben-Tal, A., A. Goryashko, A. Guslitzer, A. Nemirovski. 2004.
Adjustable robust solutions of uncertain linear programs. {\it Mathematical Programming}. {\bf 99} 351--376.

\bibitem[{Ben-Tal, Boaz and Shimrit(2009)}]{aharon2009robust}
Ben-Tal, A., G. Boaz, S. Shimrit.  2009.  Robust multi-echelon multi-period inventory control. {\it European Journal of Operational Research}. \textbf{199}(3) 922--935.


\bibitem[{Bertsekas and Shreve(1978)}]{berts}
Bertsekas, D.P., S.E. Shreve. 1978.
{\em Stochastic Optimal Control: The Discrete Time Case}.  Academic Press, New York.


\bibitem[{Bertsimas, Iancu and Parrilo(2010)}]{BIP}
Bertsimas, D., D. A. Iancu, P.A. Parrilo. 2010.
Optimality of affine policies in multistage robust optimization. {\it Mathematics of Operations Research}. {\bf 35}(2) 363-394.

\bibitem[{Bertsimas and Popescu(2005)}]{Bertsimas05}
Bertsimas, D., I. Popescu. 2005. 
Optimal inequalities in probability theory: A convex optimization approach. {\it SIAM Journal on Optimization}. {\bf 15}(3) 780-804.


\bibitem[{Bertsimas and Thiele(2006)}]{BT06}
Bertsimas, D., A. Thiele. 2006.
A robust optimization approach to inventory theory. {\it Operations Research}. \textbf{54}(1) 150-168.

\bibitem[{Boda and Filar(2006)}]{bod:06}
Boda, K., J.A. Filar. 2006.
Time consistent dynamic risk measures. {\em Mathematical Methods in Operations Research}. \textbf{63} 169--186.

\bibitem[{Bonnans and Shapiro(2000)}]{BS}
Bonnans, J. F., A. Shapiro. 2000.
{\em Perturbation Analysis of Optimization Problems}. Springer-Verlag, New York.

\bibitem[{Carpentier et al.(2012)}]{carp}
Carpentier, P.,  J.P. Chancelier, G. Cohen, M. De Lara, P. Girardeau. 2012.
Dynamic consistency for stochastic optimal control problems. {\em Annals of Operations Research}. \textbf{200} 247--263.

\bibitem[{Carrizosaa, Olivares-Nadal and Ramirez-Cobob(2016)}]{carrizosaa2014robust}
Carrizosaa, E., A. Olivares-Nadal, P. Ramirez-Cobob.  2016.
Robust newsvendor problem with autoregressive demand.  {\it Computers \& Operations Research}. \textbf{68} 123-133.

\bibitem[{Chen, Li and Guo(2013)}]{chen2013}
Chen, Z., G. Li, J. Guo. 2013.
Optimal investment policy in the time consistent mean-variance formulation. {\em Insurance: Mathematics and Economics}. \textbf{52} 145--156.

\bibitem[{Chen et al.(2007)}]{CSLS}
Chen, X., M. Sim, D. Simchi-Levi, P. Sun. 2007.
Risk aversion in inventory management. {\em Operations Research}. \textbf{55}(5) 828-842.

\bibitem[{Chen and Sim(2009)}]{ChenSim}
Chen, W., M. Sim. 2009.
Goal-driven optimization. {\it Operations Research}. {\bf 57}(2) 342-357.

\bibitem[{Chen and Sun(2012)}]{ChenSun}
Chen, X., P. Sun. 2012.
Optimal structural policies for ambiguity and risk averse inventory and pricing models. {\it SIAM Journal on Control and Optimization}. \textbf{50}(1) 133-146.

\bibitem[{Cheridito and Kupper(2009)}]{CK2009}
Cheridito, P., M. Kupper. 2009.
Recursiveness of indifference prices and translationinvariant preferences. {\em  Mathematics and Financial Economics}. \textbf{2} 173--188.


\bibitem[{Choi and Ruszczynski(2008)}]{choi2008risk}
Choi, S., A. Ruszczynski.  2008.  A risk-averse newsvendor with law invariant coherent measures of risk.
{\it Operations Research Letters}.  \textbf{36}(1) 77--82.

\bibitem[{Choi, Ruszczy\'nski and Zhao(2011)}]{CRZ:11}
Choi, S., A. Ruszczy\'nski, Y. Zhao. 2011.
A multiproduct risk-averse newswendor with law-invariant coherent measures of risk. {\em Operations Research}. {\bf 59} 346--354.

\bibitem[Delage and Iancu(2015)]{Delage15}
Delage, E., D. Iancu.  2015.
Robust Multistage Decision Making. {\em INFORMS Tutorials in Operations Research}.  {\tt http://dx.doi.org/10.1287/educ.2015.0139}.






\bibitem[{Edgeworth(1888)}]{edg:88}
Edgeworth, F. 1888.
The mathematical theory of banking. {\em Royal Statistical Society}. \textbf{51} 113--127.

\bibitem[{Epstein and Schneider(2003)}]{ES}
Epstein, L. G., M. Schneider. 2003.
Recursive multiple-priors. {\em Journal of Economic Theory}. \textbf{113}(1) 1-31.

\bibitem[{Etner, Jeleva and Tallon(2012)}]{EJT}
Etner, J., M. Jeleva, J.-M. Tallon. 2012.
Decision theory under ambiguity. {\em Journal of Economic Surveys}. \textbf{26}(2) 234-270.


\bibitem[{Federgruen and Katalan(1999)}]{Federgruen99}
Federgruen, A., Z. Katalan. 1999.
The impact of adding a maketo-order item to a make-to-stock production system. {\it Management Science}. \textbf{45}(7) 980-994.


\bibitem[{Gabrel, Murat and Thiele(2014)}]{VMT14}
Gabrel, V., C. Murat, A. Thiele.  2014.  Recent advances in robust optimization: an overview.  {\it European Journal of Operational Research}.  \textbf{235}(3) 471-483.

\bibitem[{Gallego(1998)}]{gallego1998new}
Gallego, G.  1998.  New bounds and heuristics for (Q, r) policies. {\it Management Science}. \textbf{44}(2) 219--233.

\bibitem[{Gallego(2001)}]{gallego2001minimax}
Gallego, G.  2001.  Minimax analysis for finite-horizon inventory models.  {\it IIE Transactions}.  \textbf{33}(10) 861--874.

\bibitem[{Gallego, Katircioglu and Ramachandran(2007)}]{GKR}
Gallego, G., K. Katircioglu, B. Ramachandran. 2007.
Inventory management under highly uncertain demand. {\it Operations Research Letters}. \textbf{35}(3) 281-289.

\bibitem[{Gallego and Moon(1993)}]{GM}
Gallego, G., I. Moon. 1993.
The distribution free newsboy problem: review and extensions. {\it Journal of the Operational Research Society}. \textbf{44}(8) 825--834.

\bibitem[{Gallego and Moon(1994)}]{moon1994distribution}
Gallego, G., I. Moon.  1994.  Distribution free procedures for some inventory models.  {\it Journal of the Operational research Society}.  651--658.

\bibitem[{Grunwald and Halpern(2011)}]{GH2011}
Grunwald, P. D., J. Y. Halpern. 2011.
Making decisions using sets of probabilities: updating, time consistency, and calibration. {\it Journal of Artificial Intelligence Research} \textbf{42}(1) 393-426.


\bibitem[{Hanasusanto et al.(2015)}]{hanasusanto2012distributionally}
Hanasusanto, G., A. Grani, D. Kuhn, W. Wallace, S. Zymler.  2015.
Distributionally robust multi-item newsvendor problems with multimodal demand distributions.  {\it Mathematical Programming}. \textbf{152}(1) 1-32.

\bibitem[{Hansen and Sargent(2001)}]{HS01}
Hansen, L. P., T. J. Sargent. 2001.
Robust control and model uncertainty. {\em American Economic Review}. \textbf{91}(2) 60-66.

\bibitem[{Homem-de-Mello and Pagnoncelli(2016)}]{Homem.16}
Homem-de-Mello, T., B. Pagnoncelli. 2016.
Risk aversion in multistage stochastic programming: A modeling and algorithmic perspective. {\it European Journal of Operational Research}. \textbf{249}(1) 188-199.

\bibitem[{Huang et al.(2011)}]{Huang11}
Huang, P., D. Iancu, M. Petrik, D. Subramanian. 2011.
The Price of Dynamic Inconsistency for Distortion Risk Measures.
Technical Report.


\bibitem[{Iancu, Petrik and Subramanian(2015)}]{IPS}
Iancu, D.A., M. Petrik, D. Subramanian. 2015.
Tight Approximations of Dynamic Risk Measures. {\it Mathematics of Operations Research}. \textbf{40}(3) 655-682.


\bibitem[{Ignall and Veinott(1969)}]{IV}
Ignall, E., A. Veinott. 1969.
Optimality of myopic inventory policies for several substitute products. {\it Management Science}. \textbf{15}(5) 284-304.

\bibitem[Iancu and Trichakis(2014)]{Iancu14}
Iancu, D.A., N. Trichakis. 2014.
Pareto Efficiency in Robust Optimization. {\it Management Science}. \textbf{60}(1), 130--147.


\bibitem[{Isii(1962)}]{isii}
Isii, K., 1962.
On sharpness of Tchebycheff-type inequalities. {\em Ann. Inst. Stat. Math.}, \textbf{14} 185--197.

\bibitem[{Iyengar(2005)}]{iye:05}
Iyengar, G.N.. 2005.
Robust Dynamic Programming. {\em Mathematics of Operations Research}. \textbf{30} 257--280.


\bibitem[{Jagannathan(1978)}]{Jagannathan78}
Jagannathan, R. 1978.
A minimax ordering policy for the infinite stage dynamic inventory problem. {\it Management Science}. \textbf{24} 1138-1149.


\bibitem[{Kaminsky and Kaya(2009)}]{Kaminsky09}
Kaminsky, P.,  O. Kaya. 2009. 
Combined make-to-order/make-to-stock supply chains. {\it IIE Transactions}. \textbf{41}(2) 103-119.


\bibitem[{Kasugai and Kasegai(1961)}]{kasugai1961note}
Kasugai, H., T. Kasegai.  1961.  Note on minimax regret ordering policy-static and dynamic solutions and a comparison to maximin policy.  {\it Journal of the Operations Research Society of Japan}.  \textbf{3}(4)  155--169.


\bibitem[{Klabjan, Simchi-Levi and Song(2013)}]{KSLS}
Klabjan, D., D. Simchi-Levi, M. Song. 2013.
Robust stochastic lot-sizing by means of histograms. {\it Production and Operations Management}. \textbf{22} 691-710.



\bibitem[{Lam and Ghosh(2013)}]{GL2013}
Lam, H., S. Ghosh. 2013.
Iterative methods for robust estimation under bivariate distributional uncertainty. {\it Proceedings of the Winter Simulation Conference}.

\bibitem[{Landau(1987)}]{lan:87}
Moments in mathematics, H.J. Landau (ed.), {\em Proc. Sympos. Appl. Math.}, 37, Amer. Math. Soc., Providence, RI, 1987.

\bibitem[{Levi, Perakis and Uichanco(2015)}]{LPU11}
Levi, R., G. Perakis, J. Uichanco. 2015.
The data-driven newsvendor problem: new bounds and insights. {\em Operations Research}. \textbf{63}(6) 1294-1306.

\bibitem[{Lovejoy(1992)}]{Lovejoy}
Lovejoy, W. 1992.
Stopped myopic policies for some inventory models with uncertain demand distributions. {\it Management Science}. \textbf{38}(5) 688-707.

\bibitem[{Natarajan and Zhou(2007)}]{Natarajan}
Natarajan, K., L. Zhou. 2007.
A mean-variance bound for a three-piece linear function. {\it Probability in the Engineering and Informational Sciences}. \textbf{21}(4) 611-621.


\bibitem[{Nilim and El Ghaoui(2005)}]{nil:05}
Nilim, A., L. El Ghaoui. 2005.
Robust Control of Markov Decision Processes with Uncertain Transition Matrices. {\em Operations Research}. \textbf{53} 780--798.

\bibitem[{Perakis and Roels(2008)}]{PR}
Perakis, G., G. Roels. 2008.
Regret in the newsvendor model with partial information. {\it Operations Research}. \textbf{56}(1) 188-203.

\bibitem[{Popescu(2005)}]{Popescu}
Popescu, I., 2005.
A Semidefinite Programming Approach to Optimal Moment Bounds for Convex Classes of Distributions. {\em Mathematics of Operations Research}. \textbf{50}(3) 632-657.




\bibitem[{Rajagopalan(2002)}]{Rajagopalan02}
Rajagopalan, S. 2002.
Make-to-order or make-to-stock: model and application. {\it Management Science}. \textbf{48}(2) 241-256.


\bibitem[{Riedel(2004)}]{Ri}
Riedel, F. 2004.
Dynamic coherent risk measures. {\it Stochastic Processes and their Applications}. \textbf{112} 185--200.

\bibitem[{Roorda and Schumacher(2007)}]{RS07}
Roorda, B., J.M. Schumacher. 2007.
Time consistency conditions for acceptability measures, with an applications to Tail Value at Risk. {\em Insurance: Mathematics
and Economics}. \textbf{40} 209--230.

\bibitem[{Ruszczy\'{n}ski and Shapiro(2006)}]{rs06}
Ruszczy\'{n}ski, A., A. Shapiro. 2006.  Conditional risk mappings. {\it Mathematics of Operations Research}. \textbf{31} 544--561.


\bibitem[{Ruszczy\'{n}ski(2010)}]{Ru}
Ruszczy\'{n}ski, A. 2010.
Risk-averse dynamic programming for Markov decision processes. {\it Mathematical Programming, Series B}. \textbf{125} 235--261.

\bibitem[{Scarf(1958)}]{Scarf}
Scarf, H. 1958.
A min-max solution of an inventory problem. K. Arrow, ed {\it Studies in the Mathematical Theory of Inventory and Production}. Stanford University Press,         Stanford, CA, 201--209.

\bibitem[{Scarf(1959)}]{scarf1959bayes}
Scarf, H.  1959.  Bayes solution of the statistical inventory problem.  {\it The Annals of Mathematical Statistics}.  490--508.



\bibitem[{Scarf(1960)}]{scarf1960some}
Scarf, H.  1960.  Some remarks on Bayes solutions to the inventory problem.  {\it Naval Research Logistics Quarterly}.  \textbf{7}(4)  591--596.



\bibitem[{See and Sim(2010)}]{SS}
See, C. T., M. Sim. 2010.
Robust approximation to multiperiod inventory management. {\it Operations Research}. \textbf{58}(3) 583-594.




\bibitem[{Shapiro, Dentcheva and Ruszczy\'{n}ski(2009)}]{SDR}
Shapiro, A., D. Dentcheva, A. Ruszczy\'{n}ski. 2009.
{\em Lectures on Stochastic Programming: Modeling and Theory}. SIAM, Philadelphia.


\bibitem[{Shapiro(2009)}]{S-09}
Shapiro, A.  2009.
On a time consistency concept in risk averse multistage stochastic programming. {\it Operations Research Letters}. \textbf{37} 143--147.


\bibitem[{Shapiro(2011)}]{dpadjust11}
Shapiro, A.  2011.
A dynamic programming approach to adjustable robust optimization.  {\em Operations Research Letters}. \textbf{39}(2) 83-87.

\bibitem[{Shapiro(2012)}]{S-12}
Shapiro, A.  2012.
Minimax and risk averse multistage stochastic programming. {\it European Journal of Operational Research}. \textbf{219} 719--726.

\bibitem[{Shapiro(2016)}]{sha-14}
Shapiro, A.  2016.
Rectangular sets of probability measures. {\it Operations Research}. \textbf{64}(2) 528-541.

\bibitem[{Shapiro(2017)}]{Shapiro17a}
Shapiro, A.  2017.
Interchangeability principle and dynamic equations in risk averse stochastic programming. {\it Operations Research Letters}. \textbf{45}(4) 377-381.

\bibitem[{Shapiro and Xin(2017)}]{Shapiro17}
Shapiro, A., L. Xin. 2017.
Time inconsistency of optimal policies of distributionally robust inventory models. Optimization Online 2017/09/6232.


\bibitem[{Strotz(1955)}]{Strotz55}
Strotz, R. H.. 1955.
Myopia and inconsistency in dynamic utility maximization. {\it The Review of Economic Studies}. \textbf{23}(3) 165-180.

\bibitem[{Veinott(1965)}]{Veinott65}
Veinott, Jr., A. F.. 1965.
Optimal policy for a multi-product, dynamic, nonstationary inventory problem. {\it Management Science}. \textbf{12}(3) 206-222.

\bibitem[{Wang(1999)}]{Wan99}
Wang, T. 1999.
A class of dynamic risk measure. {\em Working Paper}. University of British Columbia.

\bibitem[{Wiesemann, Kuhn and Rustem(2013)}]{WKR}
Wiesemann, W., D. Kuhn, B. Rustem. 2013.
Robust Markov decision processes. {\it Mathematics of Operations Research}. \textbf{38}(1) 153-183.


\bibitem[{Williams(1984)}]{Williams84}
Williams, T.M. 1984.
Special products and uncertainty in production/inventory systems. {\it European Journal of Operations Research}. \textbf{15} 46-54.

\bibitem[{Xin, Goldberg and Shapiro(2013)}]{Xin.13}
Xin, L., D. Goldberg, A. Shapiro. Time (in) consistency of multistage distributionally robust inventory models with moment constraints. arXiv preprint arXiv:1304.3074 (2013).


\bibitem[{Yang(2013)}]{Yang13}
Yang, J.  2013.  Inventory and Price Control under Time-consistent Coherent and Markov Risk Measure. Preprint.

\bibitem[{Yue, Chen and Wang(2006)}]{yue2006expected}
Yue, J., B. Chen, M. Wang.  2006.  Expected value of distribution information for the newsvendor problem.  {\it Operations Research}.  \textbf{54}(6) 1128--1136.



\bibitem[{Zhu, Zhang and Ye(2013)}]{zhu2013newsvendor}
Zhu, Z., J. Zhang, Y. Ye.  2013.
Newsvendor optimization with limited distribution information.  {\it Optimization methods and software}.  \textbf{28}(3) 640-667.

\bibitem[{Zipkin(2000)}]{Zi}
Zipkin, P.H. 2000.
{\em Foundations of Inventory Management}. McGraw-Hill, Boston.


\end{thebibliography}

\newpage
\section{Appendix}\label{appsec}

\subsection{Proof of Theorem\ \ref{Scarfold}}
\proof{Proof of Theorem\ \ref{Scarfold} : }
We first compute the value of $\psi(x)$ for all $x \in \bbr$, and proceed by a case analysis.  First, suppose $x < 0$.
In this case, $\mathbb{E}_{Q}[\Psi(x,D)] = c x + b (\mu - x)$ for all $Q \in \cM$, and thus
\begin{equation}\label{psineg}
\psi(x) = c x + b (\mu - x).
\end{equation}
  Now, suppose $x \geq 0$.  Then it is easily verified that
\begin{equation}\label{rewrite1}
\psi(x) = c x+\frac{(h-b)(x-\mu)}{2}+\frac{b+h}{2} \sup_{Q\in \mathfrak{M}}\mathbb{E}_{Q}\left[|x-D|\right].
\end{equation}
Hence to compute $\psi(x)$, it suffices to solve $\sup_{Q\in \mathfrak{M}}\mathbb{E}_{Q}\left[|x-D|\right]$, and we proceed by a case analysis.   Recall that $\cf(z) := \big( (z - \mu)^2 + \sigma^2 \big)^{\frac{1}{2}}$ for all $z \in \bbr$, and $\cf^{-1}(z)$ denotes the reciprocal of $\cf(z)$.
\\\\ First, suppose $x \geq \frac{\mu^2+\sigma^2}{2\mu}$.  Let us define $\bar{\lambda} = (\bar{\lambda}_0, \bar{\lambda}_1, \bar{\lambda}_2)$ such that
$$\bar{\lambda}_0 := \frac{1}{2} \big( x^2 \cf^{-1}(x) + \cf(x) \big) ,\ \ \ \bar{\lambda}_1 := - x \cf^{-1}(x),\ \ \ \bar{\lambda}_2 := \frac{1}{2} \cf^{-1}(x),$$
and let $\bar{g}(d) := \bar{\lambda}_0 + \bar{\lambda}_1 d + \bar{\lambda}_2 d^2$ for all $d \in \bbr$.
Then it follows from a straightforward calculation that $\bar{g}(d)$ and $|x-d|$ are tangent at $\bar{d}_1 := x - \cf(x)$ and $\bar{d}_2 := x + \cf(x)$, and consequently $\bar{g}(d) \geq |x - d|$ for all $d \in \bbr_+$.  Hence $\bar{\lambda}$ is feasible for the dual Problem {\rm (\ref{dual-2})}.  Also, as $x \geq \frac{\mu^2+\sigma^2}{2\mu}$ implies $\bar{d}_1 \geq 0$, it is easily verified that the probability measure $\bar{Q}$ such that
$$\bar{Q}(\bar{d}_1) = \sigma^2 \bigg( \sigma^2 + \big( x - \cf(x) - \mu \big)^2 \bigg)^{-1},\ \ \ \bar{Q}(\bar{d}_2) =
1 - \sigma^2 \bigg( \sigma^2 + \big( x - \cf(x) - \mu \big)^2 \bigg)^{-1}$$
 is feasible for the primal Problem {\rm (\ref{dual-1})}.  It follows from Proposition\ \ref{dualprop1} that $\bar{Q}$ is an optimal primal solution.  Combining the above and simplifying the relevant algebra, we conclude that in this case
\begin{equation}\label{psi1}
\psi(x) = \psi_1(x) := c\mu+\frac{b+h}{2}\cf(x)-\frac{b-h-2c}{2}(x-\mu).
\end{equation}
\ \\Alternatively, suppose $x \in [0 , \frac{\mu^2+\sigma^2}{2\mu})$.
Let us define
$\hat{\lambda} = (\hat{\lambda}_0, \hat{\lambda}_1, \hat{\lambda}_2)$ such that
$$\hat{\lambda}_0 := x\ \ \ ,\ \ \ \hat{\lambda}_1 := 1 - 4 x \mu (\mu^2 + \sigma^2)^{-1}\ \ \ ,\ \ \ \hat{\lambda}_2 := 2 x \big( \mu (\mu^2 + \sigma^2)^{-1} \big)^2,$$
and let $\hat{g}(d) := \hat{\lambda}_0 + \hat{\lambda}_1 d + \hat{\lambda}_2 d^2$ for all $d \in \bbr$.
Then it follows from a straightforward calculation that $\hat{g}(d)$ and $|x-d|$ are tangent at $\hat{d}_1 := \mu^{-1} (\mu^2 + \sigma^2)$, and intersect at $\hat{d}_2 := 0$, with $\hat{g}'(0) \geq -1$.  It follows that $\hat{g}(d) \geq |x - d|$ for all $d \in \bbr_+$.  Hence $\hat{\lambda}$ is feasible for the dual Problem {\rm (\ref{dual-2})}.  Also, it is easily verified that the probability measure $\hat{Q}$ such that
$$\hat{Q}(\hat{d}_1) = \mu^2 (\mu^2 + \sigma^2)^{-1},\ \ \ \hat{Q}(\hat{d}_2) = 1 - \mu^2 (\mu^2 + \sigma^2)^{-1}$$
is feasible for the primal Problem {\rm (\ref{dual-1})}.  It follows from Proposition\ \ref{dualprop1} that $\hat{Q}$ is an optimal primal solution.  Combining the above and simplifying the relevant algebra, we conclude that in this case
\begin{equation}\label{psi2}
\psi(x) = \psi_2(x) := \frac{(h+c)\sigma^2-(b-c)\mu^2}{\mu^2+\sigma^2}x+ b\mu.
\end{equation}
We now use the above to complete the proof of the theorem.  Note that since by assumption $b > c$, it follows from (\ref{psineg}) that $\argmin_{x \in \bbr} \psi(x) \subseteq \bbr_+$.  Recall that $\kappa = \frac{b - h - 2c}{b + h}$.  Furthermore, our assumptions, i.e. $b > c, h + c > 0$, imply that $|\kappa| < 1$.  Let $\chi := \mu + \kappa \sigma (1 - \kappa^2)^{-\frac{1}{2}}$.
It follows from a straightforward calculation that $\psi_1$ is a strictly convex function on $\bbr$, and $\psi_1(\chi) = 0$, i.e.
$\psi_1$ is strictly decreasing on $(-\infty,\chi)$, and strictly increasing on $(\chi,\infty)$.  Furthermore, it follows from a similar calculation that
\begin{equation}\label{iffry1}
\frac{\sigma^2}{\mu^2} - \frac{b-c}{h+c}\ \ \ \textrm{is the same sign as}\ \ \ \frac{\mu^2+\sigma^2}{2\mu} - \chi.
\end{equation}
We now proceed by a case analysis.  First, suppose $\frac{\sigma^2}{\mu^2}>\frac{b-c}{h+c}$.  In this case, $\psi_2$ is a linear function with strictly positive slope, and thus $\argmin_{x \in [0,\frac{\mu^2+\sigma^2}{2\mu}]} \psi(x) = \lbrace 0 \rbrace$.  Furthermore, it follows from (\ref{iffry1}) that $\chi < \frac{\mu^2+\sigma^2}{2\mu}$, which implies that $\psi_1$ is strictly increasing on $[\frac{\mu^2+\sigma^2}{2\mu},\infty)$.  It follows from the continuity of $\psi$ that
$\argmin_{x \geq \frac{\mu^2+\sigma^2}{2\mu}} \psi(x) = \lbrace \frac{\mu^2+\sigma^2}{2\mu} \rbrace$.  Combining the above, we conclude that $\argmin_{x \in \bbr} \psi(x) = \lbrace 0 \rbrace$.
\\\\Next, suppose $\frac{\sigma^2}{\mu^2}<\frac{b-c}{h+c}$.  In this case, $\psi_2$ is a linear function with strictly negative slope, and thus $\argmin_{x \in [0,\frac{\mu^2+\sigma^2}{2\mu}]} \psi(x) = \lbrace \frac{\mu^2+\sigma^2}{2\mu} \rbrace$.  Furthermore, it follows from (\ref{iffry1}) that $\chi > \frac{\mu^2+\sigma^2}{2\mu}$, which implies that
$\argmin_{x \geq \frac{\mu^2+\sigma^2}{2\mu}} \psi(x) = \lbrace \chi \rbrace$.  Combining the above, we conclude that $\argmin_{x \in \bbr} \psi(x) = \lbrace \chi \rbrace$.
\\\\Finally, suppose that $\frac{\sigma^2}{\mu^2} = \frac{b-c}{h+c}$.  In this case, $\psi_2$ is a constant function, and thus $\argmin_{x \in [0,\frac{\mu^2+\sigma^2}{2\mu}]} \psi(x) = [0,\frac{\mu^2+\sigma^2}{2\mu}]$.  Furthermore, it follows from (\ref{iffry1}) that $\chi = \frac{\mu^2+\sigma^2}{2\mu}$, which implies that
$\argmin_{x \geq \frac{\mu^2+\sigma^2}{2\mu}} \psi(x) = \lbrace \frac{\mu^2+\sigma^2}{2\mu} \rbrace$.  Combining the above, we conclude that $\argmin_{x \in \bbr} \psi(x) = [0,\frac{\mu^2+\sigma^2}{2\mu}]$.
\\\\Combining all of the above with another straightforward calculation completes the proof of the theorem. $\Halmos$
\endproof
\subsection{Proof of Proposition\ \ref{Scarfold2}}
\proof{Proof of Proposition\ \ref{Scarfold2} : }
Let $\delta := \frac{\sigma^2}{\mu^2+\sigma^2}, \tau := \frac{\mu^2+\sigma^2}{\mu}$.
Let $Q_2^*$ be the probability measure such that
\[
Q_2^*(0) = \delta ,\ \ \ Q_2^*\left(\tau\right) =
1 - \delta.
\]
Recall that $b - c > 0$, and $(h + c)\sigma^2 > (b - c)\mu^2$, which we denote by assumption A1.
Note that the value of any feasible solution $Q_1$ to Problem {\rm (\ref{ddminmax})} is at least $\mathbb{E}_{Q_1\times Q_2^*}\Big[\Psi(D_1,D_2)\Big]$, which itself equals the sum of $c\mu$ and
\begin{equation}
\mathbb{E}_{Q_1}
\Big[\bigg(
\delta \big((b-c)[0-D_1]_++(h+c)[D_1-0]_+
\big)+(1 - \delta)\big((b-c)[\tau-D_1]_+
+(h+c)[D_1 - \tau]_+\big) \bigg) I(D_1 > 0) \Big]\label{morezero}
\end{equation}
\begin{equation}
+\ \mathbb{E}_{Q_1}
\Big[\bigg(
\delta \big((b-c)[0-D_1]_++(h+c)[D_1-0]_+
\big)+(1-\delta)\big((b-c)[\tau - D_1]_+
+(h+c)[D_1 - \tau]_+\big) \bigg) I(D_1 < 0) \Big]\label{lesszero}
\end{equation}
\begin{equation}
+\ \mathbb{E}_{Q_1}
\Big[\bigg(
\delta \big((b-c)[0-D_1]_++(h+c)[D_1-0]_+
\big)+(1 - \delta) \big((b-c)[\tau - D_1]_+
+(h+c)[D_1 - \tau ]_+\big) \bigg) I(D_1 = 0) \Big]\label{iszero}
\end{equation}
Note that if $P(D_1 > 0) > 0$, then (\ref{morezero}) is at least
\begin{eqnarray}
\ &\ &\ \mathbb{E}\Big[ \frac{\sigma^2}{\mu^2+\sigma^2} (h + c) D_1 + \frac{\mu^2}{\mu^2+\sigma^2}(b-c)(\frac{\mu^2+\sigma^2}{\mu}-D_1) \big| D_1 > 0 \Big]P(D_1 > 0)\nonumber
\\&\ &\ \ \ >\ \ \ \mathbb{E}\Big[ \frac{\mu^2}{\mu^2+\sigma^2} (b -c) D_1 + \frac{\mu^2}{\mu^2+\sigma^2}(b-c)(\frac{\mu^2+\sigma^2}{\mu}-D_1) \big| D_1 > 0 \Big]P(D_1 > 0)\ \ \ \textrm{by A1}\nonumber
\\&\ &\ \ \ =\ \ \ (b - c)\mu P(D_1 > 0).\label{moreworse}
\end{eqnarray}
Similarly, if $P(D_1 < 0) > 0$, then (\ref{lesszero}) is at least
\begin{eqnarray}
\ &\ &\ \mathbb{E}\Big[ - \frac{\sigma^2}{\mu^2 + \sigma^2} (b - c) D_1 + \frac{\mu^2}{\mu^2 + \sigma^2} (b - c) (\frac{\mu^2 + \sigma^2}{\mu} - D_1) \big| D_1 < 0\Big] P(D_1 < 0)\nonumber
\\\ &\ &\ \ \ =\ \ \ \mathbb{E}\Big[ (b - c)(\mu - D_1) \big| D_1 < 0\Big] P(D_1 < 0)\ \ \ >\ \ \ (b-c)\mu P(D_1 < 0).\label{lessworse}
\end{eqnarray}
Furthermore, if $P(D_1 = 0) > 0$, then (\ref{iszero}) equals $(b-c)\mu P(D_1 = 0)$.  Combining with (\ref{moreworse}), (\ref{lessworse}), and the fact that the measure $\delta_0$ attains value $b \mu$ (by Theorem\ \ref{Scarfold}), completes the proof. $\Halmos$
\endproof

\subsection{Proof of Theorem\ \ref{useme1}}
\proof{Proof of Theorem\ \ref{useme1} : }
Recall that $\eta := \frac{1}{2} (c_1 + c_2)$, and $\cf(z) := \big( (z - \mu)^2 + \sigma^2 \big)^{\frac{1}{2}}$ for all $z \in \bbr$.  Also, letting $h_1(d) := - d + c_1, h_2(d) := d - c_2$ for all $d \in \bbr$, we have that $\zeta(d) = \max\{h_1(d), 0, h_2(d)\}$ for all $d \in \bbr$.
Let $Q$ be the probability measure described in (\ref{eq-th2}), and $\lambda = (\lambda_0,\lambda_1,\lambda_2)$ the vector described in (\ref{eq-th2-dual}).  Let $g(d) := \lambda_0 + \lambda_1 d + \lambda_2 d^2$.  We now prove that $g(d) \geq \zeta(d)$ for all $d \in \bbr$.  It follows from a straightforward calculation that $g(d)$ is tangent to $h_1(d)$ at $d_1 := \eta - \cf(\eta)$, and $g(d)$ is tangent to $h_2(d)$ at $d_2 := \eta + \cf(\eta)$.  Thus $g(d) \geq \max\big( h_1(d), h_2(d) \big)$ for all $d \in \bbr$, and to prove the desired claim it suffices to demonstrate that $g(d) \geq 0$ for all $d \geq 0$.  It is easily verified that
for all $d \in \bbr$,
\begin{equation}\label{repforg}
g(d) = \half \cf^{-1}(\eta) (d - \eta)^2 + \half \big(\cf(\eta) + c_1 - c_2 \big).
\end{equation}
Recall that
$$\frac{1}{4}(2 \mu - 3 c_1 + c_2)(3 c_2 - c_1 - 2 \mu) \leq \sigma^2,$$ which we denote by assumption A2.
It follows from another straightforward calculation that assumption A2 is equivalent to requiring that $\half \big(\cf(\eta) + c_1 - c_2 \big) \geq 0$.  Combining with (\ref{repforg}), we conclude that A2 implies $g(d) \geq 0$ for all $d \in \bbr$, completing the proof that $g(d) \geq \zeta(d)$ for all $d \in \bbr$.  Hence $\lambda$ is feasible for the dual Problem {\rm (\ref{dual-2})}.  Also, it is easily verified that $Q$  is feasible for the primal Problem {\rm (\ref{dual-1})}.  It follows from Proposition\ \ref{dualprop1} that $Q$ is an optimal primal solution, and $\lambda$ is an optimal dual solution.  That these optimal solutions are unique then follows from the second part of Proposition\ \ref{dualprop1} and a straightforward contradiction argument.  Combining the above and simplifying the relevant algebra completes the proof. $\Halmos$
\endproof

\subsection{Proof of Lemma\ \ref{lboundme}}
\proof{Proof of Lemma\ \ref{lboundme} : }
Suppose $i\in \{0,...,T\}$ and $Q_1,\ldots,Q_i$ are fixed.  As a notational convenience, for $k \in [1,T]$, let $\bbe_k[\cdot]$ denote $\bbe_{\otimes_{j=1}^{k} Q_j}[\cdot]$.  We now prove that (\ref{inductme}) holds for all $t \geq i+1$, and proceed by induction.
Our particular induction hypothesis will be that there exist $Q_{i+1},\ldots,Q_{i+n}$ such that
\begin{equation}\label{inductme22}
\bbe_{i+n}\big[\hat{\Psi}_t\big(x_t(y_t),D_t\big)\big] \geq \hat{\eta}_t\ \textrm{for all}\ t \in [i + 1,i+n].
\end{equation}
We first treat the base case $n = 1$.  It follows from Jensen's inequality, and the independence structure of the measures in $\cM$, that for any $Q_{i+1} \in \cM_{i+1}$,
$$\bbe_{i+1} \big[\hat{\Psi}_{i+1}\big(x_{i+1}(y_{i+1}),D_{i+1}\big)\big] \geq
\bbe_{Q_{i+1}} \big[\hat{\Psi}_{i+1}\big(\bbe_{i}[x_{i+1}(y_{i+1})],D_{i+1})\big].
$$
Taking $Q_{i+1}$ to be any element of $\Gamma^{\bbe_{i}[x_{i+1}(y_{i+1})]}_{i+1}$ ($\Gamma_1^{x_1(y_1)}$ if $i=0$) completes the proof for $n=1$.

Now, suppose the induction holds for some $n$.  It again follows from Jensen's inequality, and the independence structure of the measures in $\cM$, that for any $Q_{i+n+1} \in \cM_{i+n+1}$,
$$\bbe_{i+n+1} \big[\hat{\Psi}_{i+n+1}\big(x_{i+n+1}(y_{i+n+1}),D_{i+n+1}\big)\big] \geq
\bbe_{Q_{i+n+1}} \big[\hat{\Psi}_{i+n+1}\big(\bbe_{i+n}[x_{i+n+1}(y_{i+n+1})],D_{i+n+1})\big].
$$
Taking $Q_{i+n+1}$ to be any element of $\Gamma^{\bbe_{i+n}[x_{i+n+1}(y_{i+n+1})]}_{i+n+1}$
completes the induction, and the proof, where the second part of the lemma follows by letting $i = 0$.
$\Halmos$
\endproof

\subsection{Proof of Theorem\ \ref{isbsest}}
\proof{Proof of Theorem \ref{isbsest} : }
Note that under these assumptions, for any measure $Q \in {\mathcal M}$ (and in fact any non-negative joint distribution for demand), for any such base-stock policy $\pi$, w.p.1 $x^{\pi}_t(y_t) = x^*_t$ for all $t=1,...,T$.  It then follows from a straightforward induction that $\pi$ is a robust-w.p.1-optimal policy for the distributionally robust DP formulation, and furthermore for all $t =1,...,T$ and $y \leq x^*_t$,
\[
 V_t(y) = \hat{\eta}_t - c_t x^*_{t-1} + c_t D_{t-1} + \sum_{s = t+1}^T \rho^{s - t} (\hat{\eta}_s + c_s \mu_{s-1}),
 \]
and
\[
V_1(y) = \sum_{t=1}^T \rho^{t-1} \hat{\eta}_t - c_1 y + \sum_{t=1}^{T-1} \rho^t c_{t+1} \mu_t.
\]
Combining with Lemma  \ref{lboundme} and Observation \ref{givesbound} completes the proof.  $\Halmos$
\endproof

\subsection{Proof of Theorem\ \ref{consistent}}
\proof{Proof of Theorem  \ref{consistent} : }
Let $\Pi^{opt}$ denote the set of optimal policies for the multistage-static problem.  It follows from Theorem\ \ref{Scarfold}.(i) and Theorem\ \ref{isbsest} that $\Pi^0 \subseteq \Pi^{opt}$, and every policy $\pi \in \Pi^0$ is time consistent.  Thus to prove the theorem, it suffices to demonstrate that $\Pi^0 = \Pi^{opt}$, and we begin by showing that $\bar{\pi} = (\bar{x}_1,\ldots,\bar{x}_T) \in \Pi^{opt}$ implies $\bar{x}_1(y_1) = 0$.  Indeed,
it follows from Lemma\ \ref{lboundme} that $\bar{\pi} \in \Pi^{opt}$ implies
$$\sup_{Q \in \cM_1} \bbe_Q\big[ \hat{\Psi}_1\big( \bar{x}_1(y_1), D_1 \big) \big] = \hat{\eta}_1 = b_1 \mu_1.$$
That $\bar{x}_1(y_1)$ must equal 0 then follows from Theorem\ \ref{Scarfold}.

We now show that $\bar{\pi} \in \Pi^{opt}$ implies $\bar{x}_2(z) = 0$ for all $z \leq 0$.  We proceed by contradiction. Suppose that there exists $z' \leq 0$ such that $\bar{x}_{2}(z') \neq 0$.
It is easily verified that there exists $Q_1 \in \cM_1$ such that  $Q_1(-z') > 0$, and consequently for this choice of $Q_1$, $\bar{x}_{2}(y_2)$ is not a.s. equal to 0.  We conclude from Proposition\ \ref{Scarfold2} that there exists
$Q_2 \in \cM_2$ such that
$$\bbe_{Q_1 \times Q_2} \big[ \hat{\Psi}_2\big(\bar{x}_2(y_2),D_2\big)\big] > \hat{\eta}_2 = b_2 \mu_2.$$
As we have already demonstrated that $\bar{x}_1(y_1) = 0$, and $Q_1 \in \cM_1$, we conclude that
$$\bbe_{Q_1} \big[\hat{\Psi}_1\big(\bar{x}_1(y_1 , D_1\big)\big] = \hat{\eta}_1 = b_1 \mu_1.$$
Combining with Lemma\ \ref{lboundme} then yields a contradiction.  The proof that $\bar{x}_t(z) = 0$ for all $z \leq 0$ and $t \geq 3$ follows from a nearly identical argument, and we omit the details.  $\Halmos$
\endproof

\subsection{Proof of Theorem\ \ref{notweak1}}
We first characterize the set of optimal policies for the multistage-static problem.
\begin{lemma}\label{staticopt111}
$\Pi_s^{opt} = \tilde{\Pi}_s$, and the multistage-static problem has optimal value 18.
\end{lemma}
\proof{Proof : }
It follows from Observation \ref{singletonian} that $\cM_1$ consists of the single probability measure $Q_1$ such that  $Q_1(1) = Q_1(3) = \frac{1}{2}$.  Let $D_1$ denote a random variable distributed as $Q_1$.  Note that for any policy $\pi = (x_1,x_2) \in \Pi$, one has that $x_1(y_1) = x_1(10) \geq 10$.  Consequently, $\prob( x_1(y_1) \geq D_1 ) = 1$, and $|x_1(y_1) - D_1| = x_1(y_1) - D_1$ w.p.1.
It then follows from a straightforward calculation that the cost of any policy $\pi = (x_1,x_2) \in \Pi$ under the multistage-static formulation equals
\begin{equation}\label{staticopt1}
2 x_1(10) - 4 + \sup_{Q_2 \in \cM_2} \bbe_{Q_2} \Big[ \half \big ( \big| x_2\big( x_1(10) - 1 \big) - D_2 \big| + \big| x_2\big( x_1(10) - 3 \big) - D_2 \big| \big) \Big].
\end{equation}
Let $\bar{\pi} = (\bar{x}_1,\bar{x}_2)$ denote any optimal policy for the multistage-static problem, i.e. $\bar{\pi} \in \Pi_s^{opt}$.  Then it follows from
(\ref{staticopt1}) and a straightforward contradiction argument that
\begin{equation}\label{no1c}
\bar{x}_1(10) = 10.
\end{equation}
Combining (\ref{staticopt1}) and (\ref{no1c}), we conclude that
\begin{equation}\label{usethiseq1}
\big(\bar{x}_2(9),\bar{x}_2(7)\big) \in \argmin_{(x,y) : x \geq 9, y \geq 7 } \sup_{Q_2 \in \cM_2} \bbe_{Q_2} \big[ \half \big(|x - D_2| + |y - D_2| \big) \big].
\end{equation}
Furthermore, it follows from Lemma \ref{lboundme} and Theorem\ \ref{Scarfold} that
\begin{equation}
\inf_{(x,y) : x \geq 9, y \geq 7 } \sup_{Q_2 \in \cM_2} \bbe_{Q_2} \big[ \half \big(
|x - D_2| + |y - D_2|  \big) \big]\ \ \geq\ \ \sup_{Q_2 \in \cM_2} \bbe_{Q_2} \big[ |8 - D_2| \big] \ \ =\ \ 2. \label{usethiseq2}
\end{equation}
Noting that
$$ \half \big(|9 - D_2| + |7 - D_2| \big) = 1 + \max( - D_2 + 7, 0, D_2 - 9),$$
it then follows from a straightforward calculation and Theorem\ \ref{useme1} that
\begin{equation}\label{usethiseq25}
\sup_{Q_2 \in \cM_2} \bbe_{Q_2} \big[ \half \big(|9 - D_2| + |7 - D_2| \big) \big] = 2.
\end{equation}
Combining the above, we conclude that $\tilde{\Pi}_s \subseteq \Pi_s^{opt}$.  Also, it then follows from a straightforward calculation that the multistage-static problem has optimal value 18.

We now prove that $\tilde{\Pi}_s = \Pi_s^{opt}$.  Indeed, suppose for contradiction that there exists some optimal policy $\hat{\pi} = (\hat{x}_1,\hat{x}_2) \notin \tilde{\Pi}_s$.  In that case, it follows from (\ref{no1c}) and (\ref{usethiseq1}) that $\frac{1}{2}\big(\hat{x}_2(9) + \hat{x}_2(7)\big) > 8$.  However, it then follows from Jensen's inequality, Theorem\ \ref{Scarfold}, and (\ref{usethiseq2}) that
$$
\sup_{Q_2 \in \cM_2} \bbe_{Q_2} \big[ \half \big(|\hat{x}_2(9) - D_2| + |\hat{x}_2(7) - D_2| \big) \big]\ \ \geq\ \ \sup_{Q_2 \in \cM_2} \bbe_{Q_2} \big[ \big| \half \big(\hat{x}_2(9) + \hat{x}_2(7)\big) - D_2 \big| \big]
\ \ >\ \ 2.
$$
Combining with (\ref{usethiseq2}) and (\ref{usethiseq25}) yields a contradiction, completing the proof $\Halmos$
\endproof
We now (partially) characterize the set of robust-w.p.1-optimal policies for the distributionally robust DP problem.
\begin{lemma}\label{dynamopt111}
$\Pi_d^{opt} \subseteq \tilde{\Pi}_d$, and the distributionally robust DP problem has optimal value $17+\frac{\sqrt{5}}{2}$.
\end{lemma}
\proof{Proof : }
Let $\bar{\pi} = (\bar{x}_1,\bar{x}_2)$ denote any robust-w.p.1-optimal policy for the distributionally robust DP problem, i.e. $\bar{\pi} \in \Pi_d^{opt}$.  Then it again follows from a straightforward contradiction argument that
\begin{equation}\label{no1c2}
\bar{x}_1(10) = 10.
\end{equation}
It then follows from (\ref{dyn-mompol}) that
$$\bar{x}_2(9) \in \argmin_{x \geq 9} \sup_{Q_2 \in \cM_2} \bbe_{Q_2}[|x - D_2|],$$
and
$$\bar{x}_2(7) \in \argmin_{x \geq 7} \sup_{Q_2 \in \cM_2} \bbe_{Q_2}[|x - D_2|].$$
The lemma then follows from Theorem\ \ref{Scarfold} and a straightforward calculation.  $\Halmos$
\endproof

Combining Lemmas\ \ref{staticopt111}\ and\ \ref{dynamopt111} completes the proof of Theorem\ \ref{notweak1}.

\subsection{Proof of Theorem\ \ref{weaknotstrong}}
We first prove that the multistage-static problem is weakly time consistent.
\begin{lemma}\label{isweak}
The multistage-static problem is weakly time consistent, and both the multistage-static and distributionally robust DP problems have optimal value 2.
\end{lemma}
\proof{Proof : }
Note that
$$\hat{\Psi}_1(x_1,d_1) = |x_1 - d_1|,\ \ \ \hat{\Psi}_2(x_2,d_2) =  |x_2 - d_2|.$$
It follows from Observation \ref{singletonian} that $\cM_1$ consists of the single probability measure $Q_1$ such that  $Q_1(1) = Q_1(3) = \frac{1}{2}$.  It follows from Theorem\ \ref{Scarfold} and a straightforward calculation that
$$
\hat{\Gamma}_1 = [1,3]\ \ ,\ \ \hat{\Gamma}_2 = 10\ \ ,\ \ \hat{\eta}_2 = 1.
$$
Combining the above with Theorem  \ref{isbsest}, we conclude that the base-stock policy $\pi$ such that  $x_1(y) = \max\{3,y\}$, and $x_2(y) = \max\{10,y\}$ for all $y \in \bbr$, is optimal for both the multistage-static formulation and robust-w.p.1-optimal for the distributionally robust DP formulation, which have common optimal value 2. $\Halmos$
\endproof
We now prove that the multistage-static problem is not strongly time consistent.  In particular, consider the policy $\pi' = (x'_1,x'_2)$ such that
\begin{equation}\label{policypi1}
x'_1(y) = \max\{3,y\},\;\;{\rm\ \ and\ \ }\;\;
x'_2(y) =
\begin{cases}
9.9, & \text{if}\ y \leq 0,\\
\max\{10.1,y\}, & \textrm{otherwise}.
\end{cases}
\end{equation}
\begin{lemma}\label{notstrong}
The policy
$\pi' \in \Pi_s^{opt}$, but $\pi' \notin \Pi_d^{opt}$.  Consequently, the multistage-static problem is not strongly time consistent.
\end{lemma}
\proof{Proof : }
We first show that $\pi' \in \Pi_s^{opt}$.  
It follows from a straightforward calculation that the cost of $\pi'$ under the multistage-static formulation equals
\begin{equation}\label{optforstatic1b}
\bbe_{Q_1}|3 - D_1| + 0.1 + \sup_{Q_2 \in \cM_2} \bbe_{Q_2} \max\big\{9.9 - D_2, 0, D_2 - 10.1 \big\}.
\end{equation}
It is easily verified that the conditions of Theorem\ \ref{useme1} are met, and we may apply Theorem\ \ref{useme1} to conclude that
$\argmax_{Q_2 \in \cM_2} \bbe_{Q_2} \max\big\{ 9.9 - D_2, 0, D_2 - 10.1 \big\}$ is the probability measure $Q_2$ such that
$Q_2(9) =  \frac{1}{2}$, $Q_2(11) = \frac{1}{2}$.  It follows that the value of (\ref{optforstatic1b}) equals 2, and we conclude that $\pi' \in \Pi_s^{opt}$, completing the proof.

We now show that $\pi' \notin \Pi_d^{opt}$.  Suppose, for contradiction, that $\pi' \in \Pi_d^{opt}$.
It then follows from a straightforward calculation (and considering the measure $Q_1 \in \cM_1$ such that $Q_1(1) = Q_1(3) = \frac{1}{2}$) that
\begin{equation}\label{theargmin1}
9.9 \in \argmin_{x \geq 0} \sup_{Q_2 \in \cM_2} \bbe_{Q_2}[|x - D_2|].
\end{equation}
However, it follows from Theorem\ \ref{Scarfold} that the right-hand side of (\ref{theargmin1}) is the singleton $\lbrace 10 \rbrace$, completing the proof.  $\Halmos$
\endproof

Combining Lemmas\ \ref{isweak}\ and\ \ref{notstrong} completes the proof of Theorem\ \ref{weaknotstrong}.

\subsection{Proof of Theorem\ \ref{strongdiff}}
We first characterize the set of optimal policies for the multistage-static problem.
\begin{lemma}\label{staticopt111diff}
$\Pi_s^{opt} = \tilde{\Pi}$, and the multistage-static problem has optimal value 5.
\end{lemma}
\proof{Proof : }
It follows from Observation \ref{singletonian}  that $\cM_1$ consists of the single probability measure $Q_1$ such that  $Q_1(1) = Q_1(3) = \frac{1}{2}$.  In this case, the cost of any policy $\pi = (x_1,x_2) \in \Pi$ under the multistage-static formulation equals
\begin{equation}\label{staticoptb1}
 \sup_{Q_2 \in \cM_2} \bbe_{Q_2} \bigg[
\bbe_{Q_1} \Big[
2 \Big( x_2\big( x_1(0) - D_1 \big) - \big( x_1(0) - D_1 \big) \Big) + \big| x_2\big( x_1(0) - D_1 \big) - D_2 \big|
\Big] \bigg].
\end{equation}
We now prove that for any policy $\bar{\pi} = (\bar{x}_1,\bar{x}_2) \in \Pi_s^{opt}$, one has that
\begin{equation}\label{no1a}
\bar{x}_2\big( \bar{x}_1(0) - 1 \big) = \bar{x}_1(0) - 1\ \ \ \textrm{and}\ \ \
\bar{x}_2\big( \bar{x}_1(0) - 3 \big) = \bar{x}_1(0) - 3.
\end{equation}
Indeed, note that w.p.1, it follows from the triangle inequality that
\begin{eqnarray}
&\ &\ 2 \Big( x_2\big( x_1(0) - D_1 \big) - \big( x_1(0) - D_1 \big) \Big) + \big| x_2\big( x_1(0) - D_1 \big) - D_2 \big| \nonumber
\\&\ &\ \ \ =\ \ 2 \Big( x_2\big( x_1(0) - D_1 \big) - \big( x_1(0) - D_1 \big) \Big) + \big| x_2\big( x_1(0) - D_1 \big) - \big( x_1(0) - D_1) + \big(x_1(0) - D_1) - D_2 \big| \nonumber
\\&\ &\ \ \ \geq\ \ 2 \Big( x_2\big( x_1(0) - D_1 \big) - \big( x_1(0) - D_1 \big) \Big) +
\big|\big(x_1(0) - D_1) - D_2 \big| - \big| x_2\big( x_1(0) - D_1 \big) - \big( x_1(0) - D_1)\big| \nonumber
\\&\ &\ \ \ =\ \ x_2\big( x_1(0) - D_1 \big) - \big( x_1(0) - D_1 \big) +
\big| x_1(0) - D_1 - D_2 \big|. \label{endofeqn}
\end{eqnarray}
Now, suppose for contradiction that (\ref{no1a}) does not hold.  It follows that
$$\bbe_{Q_1} \big[ x_2\big( x_1(0) - D_1 \big) - \big( x_1(0) - D_1 \big) \big] > 0,$$
and combining with (\ref{endofeqn}), we conclude that (\ref{staticoptb1}) is strictly greater than
\begin{equation}\label{staticoptb2}
 \sup_{Q_2 \in \cM_2} \bbe_{Q_2} \bigg[
\bbe_{Q_1} \Big[ \big| x_1(0) - D_1 - D_2 \big| \Big] \bigg].
\end{equation}
Noting that (\ref{staticoptb2}) is the cost incurred by some policy satisfying (\ref{no1a}) completes the proof.

We now complete the proof of the lemma.  It suffices from the above to prove that
\begin{equation}\label{showargmin}
\argmin_{x_1 \in \bbr_+} \sup_{Q_2 \in \cM_2} \bbe_{Q_2}\left[ \half \big( | x_1 - 1 - D_2 | +
| x_1 - 3 - D_2 | \big) \right] = \lbrace 102 \rbrace.
\end{equation}
It follows from a straightforward calculation that as long as $x_1 \geq 3$, $(x_1 - 100)(104 - x_1) \leq 25$ and $x_1 - 2 - \big( (x_1 - 2 - 100)^2 + 25 \big)^{\frac{1}{2}}\geq 0$, which holds for all $x_1 \in [100,104]$, the conditions of Theorem\ \ref{useme1} are met.  We may thus apply Theorem \ref{useme1} to conclude that
for all $x_1 \in [100,104]$,
\begin{equation}\label{mustuse1}
 \sup_{Q_2 \in \cM_2} \bbe_{Q_2}\big [ \half \big( | x_1 - 1 - D_2 | +
| x_1 - 3 - D_2 | \big) \big]
\end{equation}
has the unique optimal solution $\hat{Q}_2$ such that
$$\hat{Q}_2\big ( x_1 - 2 - \big( (x_1 - 2 - 100)^2 + 25 \big)^{\frac{1}{2}} \big) =
25 \Big (
25 + \big (x_1 - 2 - \big( (x_1 - 2 - 100)^2 + 25 \big)^{\frac{1}{2}} - 100 \big)^2 \Big )^{-1},$$
and
$$\hat{Q}_2\big ( x_1 - 2 + \big( (x_1 - 2 - 100)^2 + 25 \big)^{\frac{1}{2}} \big ) =
1 - 25 \Big (25 + \big ( x_1 - 2 - \big( (x_1 - 2 - 100)^2 + 25 \big)^{\frac{1}{2}} - 100 \big )^2 \Big )^{-1}.$$
It then follows from a straightforward calculation that for $x_1 \in [100,104]$, (\ref{mustuse1}) has the value
$$g(x_1) := \big( x^2_1 - 204 x_1 + 10429\big)^{\frac{1}{2}}.$$
It is easily verified that $g$ is a strictly convex function on $[100,104]$, $g$ has its unique minimum on that interval at the point 102, and $g(102) = 5$.  The desired result then follows from the fact that (\ref{mustuse1}) is a convex function of $x_1$ on $\bbr$. $\Halmos$
\endproof
We now prove that the multistage-static problem is strongly time consistent.
\begin{lemma}\label{isstrongagain}
The multistage-static problem is strongly time consistent, and the optimal value of the distributionally robust DP problem equals $\sqrt{26}$.
\end{lemma}
\proof{Proof : }
First, we note that as in the multistage-static setting, any policy $\bar{\pi} = (\bar{x}_1,\bar{x}_2) \in \Pi_d^{opt}$ also satisfies (\ref{no1a}).  The proof is very similar to that used for the multistage-static case, and we omit the details.  To prove the lemma, it thus suffices to prove that
\begin{equation}\label{showargmin2}
\argmin_{x_1 \in \bbr_+} \left( \half \sup_{Q_2 \in \cM_2} \bbe_{Q_2}\big[ | x_1 - 1 - D_2 | \big]
+  \half \sup_{Q_2 \in \cM_2} \bbe_{Q_2}\big[ | x_1 - 3 - D_2 | \big] \right) = \lbrace 102 \rbrace.
\end{equation}
It is easily verified that for all $x_1 \in [100,104]$, we may apply Theorem\ \ref{Scarfold} to conclude that
\[
\begin{array}{ll}
 \sup\limits_{Q_2 \in \cM_2} \bbe_{Q_2}\big[ | x_1 - 1 - D_2 | \big] = \big( (x_1 - 101)^2 + 25 \big)^{\frac{1}{2}},\\
\sup\limits_{Q_2 \in \cM_2} \bbe_{Q_2}\big[ | x_1 - 3 - D_2 | \big] = \big( (x_1 - 103)^2 + 25 \big)^{\frac{1}{2}}.
\end{array}
\]
We conclude that for all $x_1 \in [100,104]$,
\begin{equation}\label{mustuse22}
\half \sup_{Q_2 \in \cM_2} \bbe_{Q_2}\big[ |x_1 - 1 - D_2 | \big]
+  \half \sup_{Q_2 \in \cM_2} \bbe_{Q_2}\big[ | x_1 - 3 - D_2 | \big]
\end{equation}
equals
\begin{equation}\label{convexsaver}
g(x_1) := \half \left (  \big( (x_1 - 101)^2 + 25 \big)^{\frac{1}{2}} +  \big( (x_1 - 103)^2 + 25 \big)^{\frac{1}{2}} \right).
\end{equation}
It is easily verified that $g(x)$ is a strictly convex function of $x$ on $[100,104]$, $g$ has its unique minimum on that interval at the point 102, and $g(102) = \sqrt{26}$.  The desired result then follows from the fact that (\ref{mustuse22}) is a convex function of $x_1$ on $\bbr$. $\Halmos$
\endproof

Combining Lemmas\ \ref{staticopt111diff}\ and\ \ref{isstrongagain} completes the proof of Theorem\ \ref{strongdiff}.

\subsection{Proof of Theorem\ \ref{nobase-stock}}
Let $\tilde{Q}_2$ denote the probability measure such that $\tilde{Q}_2(5) =
\tilde{Q}_2(11) = \frac{1}{2}$.  It may be easily verified that $\tilde{Q}_2 \in \cM_2$.  We begin by proving the following auxiliary lemma.
\begin{lemma}\label{interlem0}
$$\sup_{Q_1 \in \cM_1,\ Q_2 \in \cM_2} \bbe_{Q_1\times
Q_2}\Big[\left|10-D_1-D_2\right| \Big]= 3.$$
\end{lemma}
\proof{Proof : }
Note that
\[
   \bbe_{Q_1\times Q_2}\Big[\left|10-D_1-D_2\right| \Big]=
\bbe_{Q_2}\Big[\bbe_{Q_1}\left[\left|10-D_1-D_2\right|\Big|D_2\right] \Big].
\]
Let us define
$$\phi_{Q_1}(d)\stackrel{\Delta}{=} \bbe_{Q_1}\left[\left|10-D_1-D_2\right|\Big| \lbrace D_2 = d \rbrace \right],$$
and
$$
q(d)\stackrel{\Delta}{=}\frac{1}{6}\left(d-8\right)^2+\frac{3}{2} = \frac{73}{6} - \frac{8}{3} d + \frac{1}{6} d^2.
$$
As $\tilde{Q}_2 \in \cM_2$, to prove the lemma, it follows from Proposition\ \ref{dualprop1} that it suffices to demonstrate that for all $Q_1 \in \cM_1$, $q(5) = \phi_{Q_1}(5), q(11) = \phi_{Q_1}(11)$, and $q(d) \geq \phi_{Q_1}(d)$ for all $d \in \bbr$, as in this case for any $Q_1 \in \cM_1$, $\sup_{Q_2 \in \cM_2}E_{Q_2}[\phi_{Q_1}(D_2)] = E_{Q_2}[q(D_2)] = 3$.
We now prove that $q(d) \geq \phi_{Q_1}(d)$ for all $d \in \bbr$.  For any $Q_1 \in \cM_1$, since $10 - D_1\in \left[7-\epsilon,\ 9+\epsilon\right]$
w.p.1, it follows that $\phi_{Q_1}(d) = 10 - \mu_1 - d = 8 - d$ if $d \in [0, 7 - \epsilon]$, and $\phi_{Q_1}(d) = d + \mu_1 - 10 = d - 8$ if $d \in [9 + \epsilon,\infty)$.
It is easily verified that $q(d) - (8 - d) \geq 0$, and $q(d) - (d - 8) \geq 0$, for all $d \in \bbr$.
It follows that $q(d)\geq \phi_{Q_1}(d)$ for all $d\in (-\infty, 7-\epsilon] \bigcup [9+\epsilon, \infty)$. Noting that $\phi_{Q_1}(d)$ is
a convex function of $d$ on $(-\infty,\infty)$, we conclude that $\phi_{Q_1}(d) \leq \max\big( \phi_{Q_1}(7 - \epsilon), \phi_{Q_1}(9 + \epsilon) \big)$ for all $d \in [7 - \epsilon, 9 + \epsilon]$.  As it is easily verified that $\inf_{d\in \bbr} q(d) = \frac{3}{2}$, to prove that $q(d) \geq \phi_{Q_1}(d)$ for $d\in \left[7-\epsilon, 9+\epsilon\right]$, it suffices to show that $\max\big( \phi_{Q_1}(7 - \epsilon), \phi_{Q_1}(9 + \epsilon) \big) \leq \frac{3}{2}$.  As $\phi_{Q_1}(7 - \epsilon) = 8 - (7 - \epsilon) = 1 + \epsilon < \frac{3}{2}$, and $\phi_{Q_1}(9 + \epsilon) = (9 + \epsilon) - 8 = 1 + \epsilon < \frac{3}{2}$, combining the above we conclude that $q(d) \geq \phi(d)$ for all $d \in \bbr$.  As it is easily verified that $q(5) = \phi_{Q_1}(5) = 3$ and $q(11) = \phi_{Q_1}(11) = 3$, combining the above completes the proof.
\Halmos
\endproof

\proof{Proof of Theorem \ref{nobase-stock} : }
Note that the cost under any policy $\pi = (x_1,x_2) \in \Pi$ under the multistage-static formulation equals
$$
\sup_{Q_1 \in \cM_1,\ Q_2 \in \cM_2} \bbe_{Q_1\times Q_2}\Big[ 2 | x_1(y_1) - D_1 | + |x_2(D_1) - D_2| \big].$$
As $D_1 \leq 3 + \epsilon \leq 10 - \epsilon$ w.p.1, and $x_1(y_1) \geq y_1 = 10 - \epsilon$, we conclude that w.p.1
$$| x_1(y_1) - D_1 | = x_1(y_1) - D_1 \geq 10 - \epsilon - D_1.$$
Combining with the fact that $\mu_1 = 2$, we conclude that
$$\bbe_{Q_1\times Q_2}\Big[ 2 | x_1(y_1) - D_1 | \big] \geq 2 \big(10 - \epsilon - 2) = 2 (8 - \epsilon).$$
As $\frac{\sigma^2_2}{\mu^2_2} = \frac{9}{64} <\frac{b_2}{h_2} = 1$, and $\big( h_2 b_2 \big)^{\frac{1}{2}} \sigma_2 = 3$, it follows from Lemma\ \ref{lboundme} and Theorem\ \ref{Scarfold} that
$$\bbe_{Q_1\times Q_2}\Big[ | x_2(D_1) - D_2 | \big] \geq 3.$$
Combining the above, we conclude that the cost incurred under any policy $\pi$ is at least $19 - 2 \epsilon$.
\\\indent We now show that the cost incurred under any such policy $\tilde{\pi}$ achieves this bound, and is thus optimal.  In particular,
$$\sup_{Q_1 \in \cM_1,\ Q_2 \in \cM_2} \bbe_{Q_1\times Q_2}\Big[ 2 | \tilde{x}_1(y_1) - D_1 | + |\tilde{x}_2(D_1) - D_2| \big]$$
equals
\begin{eqnarray*}
\ &\ &\ \sup_{Q_1 \in \cM_1,\ Q_2 \in \cM_2} \bbe_{Q_1\times Q_2}\Big[ 2 | 10 - \epsilon - D_1 | + |10 - D_1 - D_2| \big]
\\&\ &\ \ \ \ \ \ = \sup_{Q_1 \in \cM_1,\ Q_2 \in \cM_2} \bbe_{Q_1\times Q_2}\Big[ 2 ( 10 - \epsilon - D_1 ) + |10 - D_1 - D_2| \big]
\\&\ &\ \ \ \ \ \ = 2 (10 - \epsilon - \mu_1) + \sup_{Q_1 \in \cM_1,\ Q_2 \in \cM_2} \bbe_{Q_1\times Q_2}\Big[ |10 - D_1 - D_2| \big]\ \ \ =\ \ \ 19 - 2 \epsilon,
\end{eqnarray*}
where the final equality follows from Lemma\ \ref{interlem0}.
\\\indent Next we show that there is no optimal base-stock policy, i.e. no base-stock policy belongs to $\Pi_s^{opt}$.  Indeed, let us suppose for contradiction that $\hat{\pi}$ is a
base-stock policy with constants $\hat{x}_1$, $\hat{x}_2$.  The cost incurred under such a policy $\hat{\pi}$ equals
$$
\sup_{Q_1 \in \cM_1,\ Q_2 \in \cM_2} \bbe_{Q_1\times Q_2}\Big[ 2 | \max(\hat{x}_1,y_1) - D_1 | + \big|\max\big( \max(\hat{x}_1,y_1) - D_1 , \hat{x}_2 \big) - D_2 \big| \big].$$
It follows from the fact that $D_1 \leq 3 + \epsilon < 10 - \epsilon$ w.p.1 for all $Q_1\in \cM_1$, and a straightforward contradiction argument (the details of which we omit), that
$\hat{\pi}$ cannot be optimal unless $\hat{x}_1\leq 10 - \epsilon$, in which case repeating our earlier arguments, we conclude that $\max(\hat{x}_1,y_1) = 10 - \epsilon$, and for any $Q_1 \in \cM_1, Q_2 \in \cM_2$,
$$\bbe_{Q_1\times Q_2}\big[ 2 | \max(\hat{x}_1,y_1) - D_1 | \big] = 2 (8 - \epsilon).$$
Thus to prove the desired claim, it suffices to demonstrate that
\begin{equation}\label{showthis0}
   \inf_{\hat{x}_2\in \bbr}\sup_{Q_1 \in \cM_1, Q_2 \in \cM_2} \bbe_{Q_1\times
Q_2}\Big[\Big| \max\{10-\epsilon-D_1, \hat{x}_2\} - D_2 \Big|\Big]> 3.
\end{equation}
We treat two different cases: $\hat{x}_2\in(-\infty, 7+ \frac{1}{2}\epsilon]$ and
$\hat{x}_2\in [7+ \frac{1}{2}\epsilon, \infty)$. If $\hat{x}_2\leq 7+
\frac{1}{2}\epsilon$, let the probability measure $\tilde{Q}_1$ be such that
$\tilde{Q}_1(1) = \tilde{Q}_1(3) = \frac{1}{2}$, where it is easily verified that $\tilde{Q}_1 \in \cM_1$.  In this case,
\begin{equation}\label{useme0000}
\sup_{Q_1 \in \cM_1, Q_2 \in \cM_2} \bbe_{Q_1\times Q_2}\Big[ \Big|\max\{10-\epsilon-D_1, \hat{x}_2\}-D_2\Big| \Big]
\end{equation}
is at least
\begin{eqnarray}
\ &\ &\ \sup_{Q_2 \in \cM_2} \bbe_{\tilde{Q}_1\times Q_2}\Big[\Big| \max\{10-\epsilon-D_1, \hat{x}_2\}-D_2\Big| \Big]\nonumber
\\&\ &\ \ \ = \sup_{Q_2 \in \cM_2} \bbe_{Q_2}\Big[ \frac{1}{2}\left|\max\{ 7-\epsilon, \hat{x}_2\} - D_2 \right| +  \frac{1}{2}\left| 9- \epsilon - D_2 \right| \Big], \label{tc:eq02}
\end{eqnarray}
where the final equality follows from the fact that $\hat{x}_2\leq 7+
\frac{1}{2}\epsilon$ implies $\max\{9-\epsilon, \hat{x}_2\} = 9 - \epsilon$.
It follows from convexity of the absolute value function that \eqref{tc:eq02} is at least
\begin{equation}\label{intereqq1}
   \sup_{Q_2 \in \cM_2} \bbe_{Q_2}\Big[ \left|\frac{1}{2}\max\{ 7-\epsilon,
\hat{x}_2\}+\frac{1}{2}(9- \epsilon)  - D_2 \right|\Big].
\end{equation}
Note that
\begin{eqnarray}
 \frac{1}{2}\max\{ 7-\epsilon, \hat{x}_2\}+\frac{1}{2}(9- \epsilon) &\geq& \frac{1}{2}(7 - \epsilon)+ \frac{1}{2}(9- \epsilon) \nonumber
\\&=& 8 - \epsilon. \label{lb11}
\end{eqnarray}
Letting $z \stackrel{\Delta}{=} \frac{1}{2}\max\{ 7-\epsilon,\hat{x}_2\}+\frac{1}{2}(9- \epsilon)$, note that (\ref{intereqq1}) equals
$\sup_{Q_2 \in \cM_2} \bbe_{Q_2}\Big[(z - D_2)^+ + (D_2 - z)^+\big].$  Applying Theorem\ \ref{Scarfold} with $c = 0, b = h = 1$, and noting that
$\frac{\mu_2^2+\sigma_2^2}{2\mu_2} = \frac{73}{16} < 8 - \epsilon = z$, we conclude that (\ref{intereqq1}) equals
\begin{equation}\label{prelim00}
\bigg( \big(  \frac{1}{2}\max\{ 7-\epsilon,\hat{x}_2\}+\frac{1}{2}(9- \epsilon) - 8 \big)^2 + 9 \bigg)^{\frac{1}{2}}.
\end{equation}
Combining (\ref{lb11}) with the fact that
\begin{eqnarray*}
 \frac{1}{2}\max\{ 7-\epsilon, \hat{x}_2\}+\frac{1}{2}(9- \epsilon) &\leq& \frac{1}{2}(7 + \frac{1}{2}\epsilon)+ \frac{1}{2}(9 - \epsilon)
\\&=& 8 - \frac{1}{4}\epsilon,
\end{eqnarray*}
we conclude that (\ref{prelim00}) is strictly greater than 3, completing the proof of (\ref{showthis0}) for the case $\hat{x}_2 \leq 7 + \frac{1}{2} \epsilon$.
\\\indent Alternatively, if $\hat{x}_2\geq 7+ \frac{1}{2}\epsilon$, let the probability
measure $\tilde{Q}_1$ be such that $\tilde{Q}_1(\frac{1+2\epsilon}{1+\epsilon}) =
\frac{(1+\epsilon)^2}{(1+\epsilon)^2+1}$ and $\tilde{Q}_1(3+\epsilon) =
\frac{1}{(1+\epsilon)^2+1}$.  Again, it is easily verified that $\tilde{Q}_1 \in \cM_1$.
In this case, (\ref{useme0000})
is at least
\begin{equation}\label{tc:eq03}
\begin{aligned}
   \sup_{Q_2 \in \cM_2} \bbe_{Q_2}\left[ \frac{1}{(1+\epsilon)^2+1} \left| \hat{x}_2
- D_2 \right| + \frac{(1+\epsilon)^2}{(1+\epsilon)^2+1} \left|
\max\left\{10-\epsilon-\frac{1+2\epsilon}{1+\epsilon}, \hat{x}_2\right\} - D_2
\right| \right].
\end{aligned}
\end{equation}
It follows from convexity of the absolute value function that \eqref{tc:eq03} is at least
\begin{equation}\label{tc:eq04}
\begin{aligned}
   \sup_{Q_2 \in \cM_2} \bbe_{Q_2}\left[ \left|\frac{1}{(1+\epsilon)^2+1}\hat{x}_2+ \frac{(1+\epsilon)^2}{(1+\epsilon)^2+1} \max\left\{10-\epsilon-\frac{1+2\epsilon}{1+\epsilon}, \hat{x}_2\right\} - D_2
\right|\right].
\end{aligned}
\end{equation}
Letting $z \stackrel{\Delta}{=} \frac{1}{(1+\epsilon)^2+1}\hat{x}_2+ \frac{(1+\epsilon)^2}{(1+\epsilon)^2+1} \max\left\{10-\epsilon-\frac{1+2\epsilon}{1+\epsilon}, \hat{x}_2\right\}$, note that (\ref{tc:eq04}) equals
$$\sup_{Q_2 \in \cM_2} \bbe_{Q_2}\Big[(z - D_2)^+ + (D_2 - z)^+\big].$$  Furthermore,
\begin{eqnarray}
&\ &\ \frac{1}{(1+\epsilon)^2+1}\hat{x}_2+  \frac{(1+\epsilon)^2}{(1+\epsilon)^2+1}
\max\left\{10-\epsilon-\frac{1+2\epsilon}{1+\epsilon}, \hat{x}_2\right\} \nonumber
\\&\ &\ \ \ \ \ \ \geq \frac{1}{(1+\epsilon)^2+1}\left(7+ \frac{1}{2}\epsilon\right)+
\frac{(1+\epsilon)^2}{(1+\epsilon)^2+1}
\left(10-\epsilon-\frac{1+2\epsilon}{1+\epsilon}\right) \nonumber
\\&\ &\ \ \ \ \ \ = 8+\frac{\frac{1}{2}-2\epsilon-\epsilon^2}{(1+\epsilon)^2+1}\epsilon. \label{useifneeded1}
\end{eqnarray}
Applying Theorem\ \ref{Scarfold} with $c = 0, b = h = 1$, and noting that
$\frac{\mu_2^2+\sigma_2^2}{2\mu_2} = \frac{73}{16} < 8+\frac{\frac{1}{2}-2\epsilon-\epsilon^2}{(1+\epsilon)^2+1}\epsilon = z$ (having applied (\ref{tc-epsilon})),
we conclude that (\ref{tc:eq04}) equals
\begin{equation}\label{prelim111}
\bigg( \big(   \frac{1}{(1+\epsilon)^2+1}\hat{x}_2+ \frac{(1+\epsilon)^2}{(1+\epsilon)^2+1} \max\left\{10-\epsilon-\frac{1+2\epsilon}{1+\epsilon}, \hat{x}_2\right\} - 8 \big)^2 + 9 \bigg)^{\frac{1}{2}}.
\end{equation}
Combining with (\ref{useifneeded1}) and (\ref{tc-epsilon}), we conclude that (\ref{prelim111}) is strictly greater than 3, completing the proof of (\ref{showthis0}) for the case $\hat{x}_2 \leq 7 + \frac{1}{2} \epsilon$, which completes the proof.
$\Halmos$
\endproof

\section*{Acknowledgments}
The authors sincerely thank Alex Shapiro for initiating this research, and his continued support and insights throughout this project.  Research of David A. Goldberg was partly supported by NSF award CMMI 1757394.
\end{document}